\documentclass[11pt]{article}
\usepackage{amssymb,latexsym,amsfonts}
\usepackage{amsthm}
\usepackage{amsmath}
\author{Delia Letang}
\title{Subconvexity bounds in depth-aspect for automorphic $L$-functions on $GL_2$}
\date{}

\begin{document}

\maketitle

\begin{center}{\it From a spectral identity we obtain asymptotics with error-term for the second integral moments of families of automorphic L-functions for $GL_2$ over an arbitrary number field according to twists by idele characters $\chi$ with arbitrary ramification at a fixed finite place. The power-saving in the error term breaks convexity at this non-archimedean place.}\end{center}

\footnote{This research was partially supported by NSF grant DMS-0652488.}
\noindent
{\bf 1.1 INTRODUCTION}\\

\indent
The {\it convexity} or {\it trivial bound} for the zeta function is
$$|\zeta(\frac{1}{2}+it)| \ll |t|^{\frac{1}{4}+\epsilon}$$
Any improvement over $\frac{1}{4}$ in this upper bound ``breaks convexity''. Various authors have obtained subconvexity bounds in different aspects. [Weyl 1921] gave a subconvex bound 
$$|\zeta(\frac{1}{2}+it)| \ll |t|^{\frac{1}{6}+\epsilon}$$
[Burgess 1962] broke convexity in the conductor aspect for Dirichlet $L$-functions over $\mathbb{Q}$. Subconvexity bounds were also obtained for $GL_2$ $L$-functions in [Good 1982, 1986], [Meurman 1987] and [Duke-Friedlander-Iwaniec 1993, 1994, 2001]. In recent years, subconvexity results were obtained by several authors including Kowalski, Michel, Vanderkam and Venkatesh (see [Kowalski-Michel-Vanderkam 2002] and [Michel-Venkatesh 2006]).\\
\indent
Until recently, all of these results concerned integral moments of automorphic $L$-functions over $\mathbb{Q}$, or over quadratic extensions of $\mathbb{Q}$, and not over an arbitrary number field. In 2006, Diaconu and Goldfeld [Diaconu-Goldfeld 2006a, 2006b] reconsidered the cases of groundfield $\mathbb{Q}$ or complex quadratic extensions. Then Diaconu and Garrett [Diaconu-Garrett 2008] obtained asymptotics with error-term for second integral moments of $GL_2$ automorphic $L$-functions over an arbitrary number field, by a spectral  identity. In the relevant spectral identity, Diaconu-Garrett obtained asymptotics with power-saving in the error term for averages not only on the critical line but also over families of twists by gr\"ossencharakters: 
$$\sum_\chi \int_{-\infty}^\infty |L(\frac{1}{2}+it,f \otimes \chi)|^{2}\,M_\chi(t)\,dt$$
where $M_\chi(t)$ are smooth weights. They showed that this breaks convexity in the $t$-aspect.

Here we take Diaconu-Garrett's ideas in a different direction. Fixing a $GL_2$ automorphic $L$-function over a number field, we arbitrarily deform the data associated with a fixed {\it non-archimedean} place $v_1$, and allow $\chi$ to have arbitrary ramification at $v_1$.  Thus, the weights in the moment expansion are obtained from the archimedean data as well as the data associated with the ramification at the finite prime $v_1$. We then obtain asymptotics for that second moment expansion and break convexity in the $\chi$-depth-aspect at the non-archimedean place $v_1$.\\

\vspace{4mm}
\noindent
{\bf 2.1 THE MAIN RESULT}\\

\indent
In this paper we break convexity in the $\chi$-depth-aspect for a family of $L$-functions $L(\frac{1}{2}+it, f \otimes \chi)$, where $\chi$ has arbitrary ramification at a fixed finite prime $v_1$. For a cuspform $f$ on $GL_2(k)$, where $k$ is a number field of degree $d$ over $\mathbb{Q}$, the $\chi$-depth-aspect convexity bound for the twisted $L$-function $L(\frac{1}{2}+it,f \otimes \chi)$ is
$$L(\frac{1}{2}+it,f \otimes \chi) \ll q^{N(\frac{d}{2}+\epsilon)}$$
where $q^N$, with $N \geq 1$, is the conductor of $\chi$. We break convexity by decreasing the exponent, proving that
$$L(\frac{1}{2}+it,f \otimes \chi)\ll (q^N)^{\frac{d-1+\vartheta}{2}+\epsilon}$$
for $\vartheta <1$.

Subconvexity results have been obtained in different aspects by some authors. In particular, for standard $L$-functions for $GL_2$ over a number field, Diaconu and Garrett used automorphic spectral theory to break convexity in the $t$-aspect; in a recent preprint, Michel and Venkatesh have claimed a joint convexity bound simultaneously in $t$-aspect, conductor and spectral aspects by using highly non-trivial methods including ergodic theory and regularization of integrals of automorphic forms. Like Diaconu-Garrett, we use a more conceptual argument of spectral identities to break convexity in the $\chi$-depth-aspect at one place. Further research could involve breaking convexity in the depth-aspect at more than one place, and this would involve more careful consideration of analytical details.

A subconvex bound is not merely an improvement of an exponent but also has some important applications. A very concrete example is the sums-of-three-integer-squares problem. An even more striking feature is that many useful corollaries of the Grand Riemann and Lindel\"of Hypotheses are also implied by subconvexity results. Thus subconvexity bounds are sufficient for providing solutions to some natural, yet apparently unrelated, questions.\\

\vspace{4mm}
\noindent
{\bf 3.1 THE MOMENT EXPANSION}\\
\noindent
{\bf 3.1.1 Prologue}\\

\indent
In this section, the integral moment expansion is obtained by unwinding the integral representation
$$\int_{Z_\mathbb{A} G_k \backslash G_\mathbb{A}} P \acute{e} \cdot |f|^2 \,$$
where $P \acute{e}$ is a Poincar\'e series and $f$ is a cuspform on $GL_2$. We use \mbox{Diaconu} and \mbox{Garrett's} ideas (see section 2 in [Diaconu-Garrett 2008]) to reformulate the Poincar\'e series as a single object. The moment expansion  is a sum of weighted integrals of $L$-functions $L(s,f \otimes \chi)$ of twists of $f$ by idele class characters $\chi$. The weight functions depend on archimedean data and data associated with the finite place $v_1$  where $\chi$ has arbitrary ramification. We will then obtain asymptotics from the weight functions.\\

\vspace{4mm}
\noindent
{\bf 3.1.2 Unwinding to an Euler Product}\\

\indent
Define the following subgroups of $G=GL_2$:
$$P=\bigl\{ (\begin{smallmatrix} * & *\\0 & * \end{smallmatrix})\bigr \},\,\, N=\bigl\{ (\begin{smallmatrix} 1 & *\\0 & 1 \end{smallmatrix})\bigr \}, \,\, H=\bigl\{ (\begin{smallmatrix} * & 0\\0 & 1 \end{smallmatrix})\bigr \},\,\, Z=\,\,\mbox{center of}\,\,G, \,\, M=ZH=\bigl\{ (\begin{smallmatrix} * & 0\\0 & * \end{smallmatrix})\bigr \}$$
For any place $v$ of $k$, let $K_v^{\max}$ be the standard maximal compact subgroup. So for finite $v$,
$$K_v^{\max}=GL_2(\mathfrak{o}_v)$$ 
and for infinite $v$,
\begin{equation*}
K_v^{\max}=
\begin{cases}
O_2 & \text{($v \approx \mathbb{R}$)}\\
U_2 & \text{($v \approx \mathbb{C}$)}
\end{cases}
\end{equation*}

The Poincar\'e series $P \acute{e}$ is of the form
$$P \acute{e}(g)=\sum_{\gamma \in M_k \backslash G_k} \varphi(\gamma g)\,\,\,\,\,\,\,\,\,\,(\mbox{where}\,\,g \in G_\mathbb{A})$$
for suitable functions $\varphi$ on $G_{\mathbb{A}}$ defined as follows. Let 
$$\varphi=\otimes_v\, \varphi_v$$
where for finite primes $v \neq v_1$,
\begin{align*}
\varphi_v (g)=
\begin{cases}
\chi_{0,v} (m)=\left|\frac{a}{d}\right|_v^{s^\prime} & \text{(for $g=mk,\,\, m= (\begin{smallmatrix} a & 0\\0 & d \end{smallmatrix}) \in M, \,\, s^\prime \in \mathbb{C},\,\, k \in K_v^{\max}$)}\\
0 & \text{(otherwise)}
\end{cases}
\end{align*}
For finite $v=v_1$ (at which $\chi$ is allowed to be ramified)
$$\varphi_v (mg)= \left|\frac{a}{d}\right|_v^{s^\prime}\,\cdot \,\varphi_v(g)\,\,\,\,\,\,\, (m \in M_v,\, g \in G_v)$$
The data determining $\varphi_v$ for $v=v_1$ consists of its values on $N_v$ where our simple choice is
\begin{equation*}
\varphi_v \bigl(\begin{smallmatrix}1 & x\\ 0 & 1 \end{smallmatrix}\bigr)=
\begin{cases}
1 \,\,\,\,\,\,\,\,\, \hfill \text{(for $x \in \mathfrak{o}_v$)}\\
|x|_v^{-w^\prime} \,\,\,\,\,\,\,\,\, \text{(for $w^\prime \in \mathbb{C},\,\,x \not \in \mathfrak{o}_v$)}
\end{cases}
\end{equation*}
For infinite $v$ require right $K_v$-invariance and left equivariance:
$$\varphi_v (mg)= \left|\frac{a}{d}\right|_v^{s^\prime}\,\cdot \,\varphi_v(g)\,\,\, \,\,\,\,(m \in M_v,\, g \in G_v)$$ 
where
\begin{equation*}
\varphi_v \bigl(\begin{smallmatrix}1 & x\\ 0 & 1 \end{smallmatrix}\bigr)=
\begin{cases}
(1+|x|^2)^{-\frac{w}{2}}\,\,\,\,\,\, \hfill \text{(for $v \approx \mathbb{R},\,\, w \in \mathbb{C}$)}\\
(1+x\overline{x})^{-w} \,\,\,\,\,\,\, \text{(for $v \approx \mathbb{C}$)}
\end{cases}
\end{equation*}
The Poincar\'e series $P \acute{e}$ converges absolutely and locally uniformly for \mbox{$\Re(s^\prime)>1$,} \mbox{$\Re(w)>1$ for all $v|\infty$,} and for $\Re(w^\prime)>1$ (see Proposition 2.6 in [Diaconu-Garrett 2008]).\\

\noindent
We want to show that 
$$\int_{Z_\mathbb{A} G_k \backslash G_\mathbb{A}} P \acute{e} \cdot |f|^2 \, dg$$
is an integral of products of local factors of standard $L$-functions. First, the Fourier expansion of a cuspform $f$ on $G_\mathbb{A}$ is
$$f(g)=\sum_{\xi \in Z_k \backslash M_k} W_f(\xi g)$$
where $W_f$ is the Whittaker function of $f$ and $W_f=\otimes_v\, W_{f,v}$ is the factorization of $W_f$ into local data. So\\

$\displaystyle\int_{Z_\mathbb{A} G_k \backslash G_\mathbb{A}} P \acute{e} \cdot |f|^2 \, dg\,= \,\displaystyle\int_{Z_\mathbb{A} G_k \backslash G_\mathbb{A}} \displaystyle\sum_{\gamma \in M_k \backslash G_k} \varphi(\gamma g)\, |f(g)|^2 \, dg\,= \,\displaystyle\int_{Z_\mathbb{A} M_k \backslash G_\mathbb{A}} \varphi(g)\,|f(g)|^2 \, dg$\\

$= \displaystyle\int_{Z_\mathbb{A} M_k \backslash G_\mathbb{A}} \varphi(g)\,\displaystyle\sum_{\xi \in Z_k \backslash M_k} W_f(\xi g)\,\overline{f}(g) \, dg\,=\, \displaystyle\int_{Z_\mathbb{A}\backslash G_\mathbb{A}} \varphi(g)\,W_f(g)\,\overline{f}(g) \,dg$\\

\noindent
Let $C$ be the idele class group $GL_1(k) \backslash GL_1(\mathbb{A})$ and $\hat{C}$ its dual. \mbox{$\hat{C} \approx \mathbb{R} \times \hat{C}_0$} where $\hat{C}_0$ is discrete. The Mellin transform and inversion are
\begin{align*}
f(x)&=\int_{\hat{C}}\,\int_C f(y)\,\chi^{-1}(y)\,dy\,\chi(x)\,d\chi\\
&=\sum_{\chi^\prime \in \hat{C}_0} \frac{1}{2 \pi i} \int_{\Re(s)=\sigma} \int_C f(y)\,{\chi^\prime}^{-1}(y)\,|y|^{-s}\,dy\,\chi^\prime(x)\,|x|^s\,ds
\end{align*}

\noindent
With $Z_\mathbb{A} M_k \backslash M_\mathbb{A} \approx C$, and for finite $v \neq v_1$, \\

$\displaystyle\int_{Z_\mathbb{A}\backslash G_\mathbb{A}} \varphi(g)\,W_f(g)\,\overline{f}(g) \,dg\,= \,\displaystyle\int_{Z_\mathbb{A} \backslash G_\mathbb{A}} \varphi(g)\,W_f(g) \Bigl(\displaystyle\int_{\hat{C}}\,\displaystyle\int_{Z_\mathbb{A} M_k \backslash M_\mathbb{A}} \overline{f}(m^\prime g)\,\chi(m^\prime)\,dm^\prime\,d\chi \Bigr)\,dg$\\

$=\displaystyle\int_{\hat{C}}\, \Bigl(\displaystyle\int_{Z_\mathbb{A} \backslash G_\mathbb{A}} \varphi(g)\,W_f(g)\,\displaystyle\int_{Z_\mathbb{A} M_k \backslash M_\mathbb{A}} \displaystyle\sum_{\xi \in Z_k \backslash M_k}\overline{W}_f(\xi m^\prime g)\,\chi(m^\prime)\,dm^\prime\,dg \Bigr)\,d\chi$\\

$=\displaystyle\int_{\hat{C}}\, \Bigl(\displaystyle\int_{Z_\mathbb{A} \backslash G_\mathbb{A}} \varphi(g)\,W_f(g)\,\displaystyle\int_{Z_\mathbb{A} \backslash M_\mathbb{A}} \overline{W}_f(m^\prime g)\,\chi(m^\prime)\,dm^\prime \,dg \Bigr)\,d\chi$\\

$=\displaystyle\int_{\hat{C}} \displaystyle\prod_v \Bigl(\displaystyle\int_{Z_v \backslash G_v}\, \displaystyle\int_{Z_v \backslash M_v} \varphi_v(g_v)\,W_{f,v}(g_v)\,\overline{W}_{f,v}(m_v^\prime g_v)\,\chi_v(m_v^\prime)\,dm_v^\prime \,dg_v \Bigr)\,d\chi$\\

\noindent
Suppress finite $v \neq v_1$ and write the $v^{th}$ local integral as
$$\int_{Z \backslash G}\, \int_{Z \backslash M} \varphi(g)\,W_f(g)\,\overline{W}_f(m^\prime g)\,\chi(m^\prime)\,dm^\prime \,dg$$
Invoke the $v$-adic Iwasawa decomposition $G=MNK$ and rewrite the integral as
$$\int_{Z \backslash MNK}\, \int_{Z \backslash M} \varphi(mnk)\,W_f(mnk)\,\overline{W}_f(m^\prime mnk)\,\chi(m^\prime)\,dm^\prime \,dm\,dn\,dk$$
For simplicity, take $\varphi$ and $f$ to be right $K_v^{\max}$-invariant for finite $v \neq v_1$. This gives 
$$\int_{Z \backslash MN}\, \int_{Z \backslash M} \varphi(mn)\,W_f(mn)\,\overline{W}_f(m^\prime mn)\,\chi(m^\prime)\,dm^\prime \,dm\,dn$$
Replace $m^\prime$ by $m^\prime m^{-1}$ to get
$$\int_{Z \backslash MN}\,\int_{Z \backslash M} \varphi(mn)\,W_f(mn)\,\overline{W}_f(m^\prime n)\,\chi(m^\prime)\,\chi^{-1}(m)\,dm^\prime \,dm\,dn$$
The Whittaker function has the equivariance
$$W_f(ng)=\psi(n)\,W_f(g)\,\,\,(n \in N_\mathbb{A})$$
Thus,
$$W_f(mn)=W_f(mnm^{-1}m)=\psi(mnm^{-1})\,W_f(m)\,\,\,\, (\mbox{since}\,\, mnm^{-1} \in N)$$
and
$$\overline{W}_f(m^\prime n)=\overline{W}_f(m^\prime nm^{-1}m)=\overline{\psi}(m^\prime n{m^\prime}^{-1})\,\overline{W}_f(m^\prime)$$
so obtaining
$$\int_{Z \backslash MN}\, \int_{Z \backslash M} \varphi(mn)\,W_f(m)\,\overline{W}_f(m^\prime)\,\chi(m^\prime)\,\chi^{-1}(m)\,\psi(mnm^{-1})\,\overline{\psi}(m^\prime n {m^\prime}^{-1})\,dm^\prime \,dm\,dn$$ 
Let
$$X(m,m^\prime)=\int_N \varphi(n)\,\psi(mnm^{-1})\,\overline{\psi}(m^\prime n {m^\prime}^{-1})\,dn$$
We get
$$\int_{Z \backslash M}\, \int_{Z \backslash M} \chi_0(m)\,W_f(m)\,\overline{W}_f(m^\prime)\,\chi(m^\prime)\,\chi^{-1}(m)\,X(m,m^\prime)\,dm^\prime \,dm$$
Now $$W_f(mn)=\psi(mnm^{-1})\cdot W_f(m)$$ 
and
$$W_f(mn)=W_f(m) \cdot 1$$
by the right $K$-invariance of $W_f$. So for $W_f(m) \not = 0, \,\psi(mnm^{-1})=1$, and $X(m,m^\prime)=1$ for $m,m^\prime$ in the support of $W_f$. So \\

$\displaystyle\int_{Z \backslash M}(\chi_0 \cdot \chi^{-1})(m)\,W_f(m)\,dm\,\cdot \,\displaystyle\int_{Z \backslash M} \chi(m^\prime)\,\overline{W}_f(m^\prime)\,dm^\prime$\\

\mbox{$= L_v(\chi_{0,v} \cdot \chi_v^{-1}\,|y|_v^{\frac{1}{2}},\,f) \cdot L_v(\chi_v \,|y^\prime|_v^{\frac{1}{2}},\,\overline{f})\,\,\,\,\,\,(\mbox{where}\,\,\,\,\,m=(\begin{smallmatrix} y & 0\\0 & 1 \end{smallmatrix}),\,\,\,m^\prime=(\begin{smallmatrix} y^\prime & 0\\0 & 1 \end{smallmatrix}))$}\\

\noindent
is a product of local factors of $L$-functions at finite primes $v \neq v_1$. Thus the integral can be written as
$$I(\chi_0)=\sum_{\chi \in \hat{C}_0} \frac{1}{2 \pi i} \int_{\Re(s)=\sigma}L(\chi_0 \cdot \chi^{-1}\,|y|^{1-s},\,f)\, \cdot \,L(\chi \,|y^\prime|^s,\,\overline{f})\, \cdot\,\mathcal{K} _{v_1}(w^\prime, \chi_{v_1})\,\cdot \,\mathcal{K}_{\infty}(s,\chi_0,\chi)\,ds$$
where
$$\mathcal{K}_{\infty}(s,\chi_0,\chi)=\prod_{v|\infty} \mathcal{K}_v (s,\chi_{0,v},\chi_v)$$
and\\

$\mathcal{K}_v(s,\chi_{0,v},\chi_v)\,=\,\displaystyle\int_{Z_v \backslash M_v N_v}\, \displaystyle\int_{Z_v \backslash M_v} \varphi_v(m_v n_v)\,W_{f,v}(m_v n_v)\,\overline{W}_{f,v}(m_v^\prime n_v)\, \cdot$\\
$$\chi_v(m_v^\prime)\,|m_v^\prime|_v^{s-\frac{1}{2}}\,\chi_v^{-1}(m_v)\,|m_v|_v^{\frac{1}{2}-s}\,dm_v^\prime \,dm_v\,dn_v$$

$\mathcal{K}_{v_1}(w^\prime, \chi_{v_1})=$\\
$$\int_{k_v^\times} \, \int_{k_v^\times} \chi(y)\, |y|^s_{v_1}\,\chi^{-1}(y^\prime)\, |y^\prime|^{1-s}_{v_1} \, W\bigl(\begin{smallmatrix}y & 0\\ 0 & 1 \end{smallmatrix}\bigr)\, \overline{W}\bigl(\begin{smallmatrix}y^\prime & 0\\ 0 & 1 \end{smallmatrix}\bigr)\, \cdot \, \int_{k_v}\, \overline{\psi}(x \cdot (y-y^\prime))\,\varphi_{v_1}\bigl(\begin{smallmatrix}1 & x\\ 0 & 1 \end{smallmatrix}\bigr)\,dx \,\, dy \, dy^\prime $$
where there is arbitrary ramification of $\chi$ at finite $v=v_1$. The non-decoupled integrals $\mathcal{K}_v(s,\chi_{0,v},\chi_v)$ and $\mathcal{K}_{v_1}(w^\prime, \chi_{v_1})$, which represent the weight functions, will be subsequently computed.\\

\noindent
The Poincar\'e series $P \acute{e}$ has meromorphic continuation to a region in $\mathbb{C}^2$ containing $s^\prime=0$ and $w^\prime=1$. As a function of $w^\prime$, for $s^\prime=0$, it is holomorphic in the half-plane $\Re(w^\prime)>\frac{11}{18}$ ([Kim-Shahidi 2002], [Kim 2005]), except for $w^\prime=1$ where it has a pole of order $1$. This can be seen from the spectral decomposition of $P \acute{e}$ (in Chapter 3), and the argument presented by Diaconu-Garrett in the proof of Theorem 4.17 in [Diaconu-Garrett 2008].\\

\noindent
For $\Re(s^\prime)$ and $\Re(w^\prime)$ sufficiently large, the integral $I(\chi_0)=I(s^\prime,w^\prime)$ is
$$I(s^\prime,w^\prime)=\sum_{\chi \in \hat{C}_{0,S}} \frac{1}{2 \pi i} \int_{\Re(s)=\sigma}L(\chi^{-1}\,|.|^{s^\prime+1-s},\,f)\, \cdot\, L(\chi \,|.|^s,\,\overline{f})\, \cdot\,\mathcal{K} _{v_1}(w^\prime, \chi_{v_1})\,\cdot \,\mathcal{K}_{\infty}(s,s^\prime,w,\chi)\,ds$$
where $S$ is a finite set of places including archimedean places, and the sum is over the set $\hat{C}_{0,S}$ of characters ramified at the finite place $v_1$. $I(s^\prime,w^\prime)$ has meromorphic continuation to a region in $\mathbb{C}^2$ containing the point \mbox{$s^\prime=0,\,w^\prime=1$,} and $I(0,w^\prime)$ is holomorphic for $\Re(w^\prime)>\frac{11}{18}$ except for $w^\prime=1$ where it has a pole of order $1$.\\

\noindent
We will find asymptotics for $\mathcal{K} _{v_1}(w^\prime, \chi_{v_1})$ and $\mathcal{K}_{\infty}(s,s^\prime,w,\chi)$ in Section 2.4, shift the line of integration to $\Re(s)=\frac{1}{2}$ and set $s^\prime=0$. Thus for $\Re(w^\prime)$ sufficiently large
\begin{align*}
I(0,w^\prime)&=\sum_{\chi \in \hat{C}_0} \frac{1}{2 \pi i} \int_{-\infty}^\infty L(\chi^{-1}\,|.|^{\frac{1}{2}-it},\,f)\, \cdot\, L(\chi \,|.|^{\frac{1}{2}+it},\,\overline{f})\, \cdot\,\mathcal{K} _{v_1}(w^\prime, \chi_{v_1})\,\cdot \,\mathcal{K}_{\infty}(\frac{1}{2}+it,0,w,\chi)\,dt\\
&=\sum_\chi \frac{1}{2 \pi i} \int_{-\infty}^\infty |L(\frac{1}{2}+it,\,f \otimes \chi)|^2\, \cdot \,\mathcal{K} _{v_1}(w^\prime, \chi_{v_1})\,\cdot \,\mathcal{K}_{\infty}(\frac{1}{2}+it,0,w,\chi) \,dt
\end{align*}

\vspace{4mm}
\noindent
{\bf 3.1.3 The non-decoupled integrals}\\

\indent
The nondecoupled integral $\mathcal{K}_{v_1}(w^\prime, \chi_{v_1})$ is
$$\int_{k_v^\times} \, \int_{k_v^\times} \chi(y)\, |y|^s_v\,\chi^{-1}(y^\prime)\, |y^\prime|^{1-s}_v \, W\bigl(\begin{smallmatrix}y & 0\\ 0 & 1 \end{smallmatrix}\bigr)\, \overline{W}\bigl(\begin{smallmatrix}y^\prime & 0\\ 0 & 1 \end{smallmatrix}\bigr)\, \cdot \, \int_{k_v}\, \overline{\psi}(x \cdot (y-y^\prime))\,\varphi_v\bigl(\begin{smallmatrix}1 & x\\ 0 & 1 \end{smallmatrix}\bigr)\,dx \,\, dy \, dy^\prime $$
where there is arbitrary ramification of $\chi$ at finite $v=v_1$.

For finite $v=v_1$, define
\begin{equation*}
\varphi_v \bigl(\begin{smallmatrix}1 & x\\ 0 & 1 \end{smallmatrix}\bigr)=
\begin{cases}
1 \,\,\,\,\,\,\,\,\, \hfill \text{(for $x \in \mathfrak{o}_v$)}\\
|x|_v^{-w^\prime} \,\,\,\,\,\,\,\,\, \text{(for $x \not \in \mathfrak{o}_v$)}
\end{cases}
\end{equation*}

\noindent
Henceforth, we will suppress the $v$ for ease of notation. $\psi$ is the standard additive character which is trivial on the local integers $\mathfrak{o}$ and nontrivial on $\varpi^{-1}\mathfrak{o}$. $\chi$ is a ramified multiplicative character, i.e. $\chi$ is non-trivial on $\mathfrak{o}^\times$. $W$ is a Whittaker function which is invariant on $\mathfrak{o}^\times$ since it is spherical. The spherical Whittaker function is of the form
\begin{equation*}
W\bigl(\begin{smallmatrix}y & 0\\ 0 & 1 \end{smallmatrix}\bigr)=\\
\begin{cases}
\frac{\alpha^{n+1}-\beta^{n+1}}{\alpha-\beta},  & \text{$n \geq 0$}\\
0,  &\text{otherwise}
\end{cases}
\end{equation*}
\newcommand{\ord}{ord}
where $\alpha,\beta$ are Satake parameters and $\ord(y)=n$.\\
We will first compute the integral in $y$ and $y^\prime$, and then compute the integral in $x$. Now
$$\overline{\psi}(x(y-y^\prime))=\overline{\psi}(xy-xy^\prime)=\overline{\psi}(xy)\,\cdot\,\psi(xy^\prime)$$
Thus the integrals in $y$ and $y^\prime$ are as follows:
$$\int_{k^\times} \overline{\psi}(xy)\,\chi(y)\, |y|^s\,\, W\bigl(\begin{smallmatrix}y & 0\\ 0 & 1 \end{smallmatrix}\bigr)\,dy\,\cdot \, \int_{k^\times} \psi(xy^\prime)\, \chi^{-1}(y^\prime)\, |y^\prime|^{1-s}\, \overline{W}\bigl(\begin{smallmatrix}y^\prime & 0\\ 0 & 1 \end{smallmatrix}\bigr)\, dy^\prime$$

\noindent
Consider the integral in $y$:
$$\int_{k^\times} \overline{\psi}(xy)\,\chi(y)\, |y|^s\,\, W\bigl(\begin{smallmatrix}y & 0\\ 0 & 1 \end{smallmatrix}\bigr)\,dy$$
Let $\eta \in \mathfrak{o}^\times$. Replace $y$ with $y\eta$ to get
$$\int_{k^\times} \Bigl(\int_{\mathfrak{o}^\times}\overline{\psi}(xy\eta)\,\chi(y\eta)\,d\eta \Bigr)\, |y|^s\,\, W\bigl(\begin{smallmatrix}y & 0\\ 0 & 1 \end{smallmatrix}\bigr)\,dy$$

\noindent
Consider the inner integral:
$$\int_{\mathfrak{o}^\times}\overline{\psi}(xy\eta)\,\chi(y\eta)\,d\eta$$
Recall that $\chi$ is a ramified character. Let $N$ be the conductor of $\chi$. So $\chi$ is trivial on some subgroup $1+\mathfrak{m}^{N}$ of $k^\times$ and non-trivial on $1+\mathfrak{m}^{N-1}$ where $N\geq 1$ is the smallest such integer.\\

\noindent
As in [Weil 1974], a standard computation shows that
$$\int_{k^\times} {\psi}(xy)\,\chi(y)\,dy\,=\,0\,\,\,\,\mbox{unless}\,\,\,\ord(x)=-N$$
So we claim that our inner integral is zero unless $\ord(xy)=-N$. So \mbox{$\ord(y)=-\ord(x)-N$}. The integral 
$$\int_{\mathfrak{o}^\times}\overline{\psi}(xy\eta)\,\chi(y\eta)\,d\eta$$
is a Gauss sum. A Gauss sum where $\chi$ is a ramified multiplicative character with conductor $N$, is evaluated as follows. Let
$$\mathfrak{g}(\chi,\psi)=\int_{\mathfrak{o}^\times}\chi(x)\,\cdot\,\overline{\psi}(\frac{x}{\varpi^{N}})\,dx$$
\noindent
Then\\

$|\mathfrak{g}(\chi,\psi)|^2\,=\,\Bigl( \displaystyle\int_{\mathfrak{o}^\times} {\chi}(x)\, \overline{\psi}(\frac{x}{\varpi^{N}})\, dx \Bigr)\, \cdot \, \Bigl( \displaystyle\int_{\mathfrak{o}^\times} \overline{\chi}(y)\, \psi(\frac{y}{\varpi^{N}})\,dy \Bigr)$\\

$=\displaystyle\int_{\mathfrak{o}^\times}\,\displaystyle\int_{\mathfrak{o^\times}} \chi(xy^{-1})\, \psi(\frac{y-x}{\varpi^{N}})\,dx\,dy\,=\,\displaystyle\int_{\mathfrak{o}^\times}\,\displaystyle\int_{\mathfrak{o^\times}} \chi(x)\, \psi(\frac{y(1-x)}{\varpi^{N}})\,dx\,dy\,\,\,\,\mbox{by replacing $x$ with $xy$}$\\

$=\displaystyle\int_{\mathfrak{o}^\times}\chi(x)\,\displaystyle\int_\mathfrak{o} \psi(\frac{y(1-x)}{\varpi^{N}})\,dx\,dy\, - \,\displaystyle\int_{\mathfrak{o}^\times}\, \displaystyle\int_\mathfrak{m}\chi(x)\,\psi(\frac{y(1-x)}{\varpi^{N}})\,dx\,dy$\\ 

\noindent
In the first integral, since $\chi$ is trivial on $1+\varpi^{N}\mathfrak{o}$, only $x \in 1+\varpi^{N} \mathfrak{o}$ will contribute, and since $\psi$ is trivial on $\mathfrak{o}$, the first integral is
$$\mu(1+\varpi^{N} \mathfrak{o}) \cdot \mu(\mathfrak{o})= \frac{q^{2-N}}{(q-1)^2}$$
Note that
$$\mu(1+\varpi^{N} \mathfrak{o})\,=\, \frac{\mu(\mathfrak{o}^\times)}{[\mathfrak{o}^\times:1+\varpi^{N} \mathfrak{o}]}\,=\, \frac{\mu(\mathfrak{o}^\times)}{(q-1)q^{N-1}}\,=\, \frac{1}{(q-1)q^{N-1}}$$
where $\mu(\mathfrak{o^\times})$ is normalized to $1$ and $\mu(\mathfrak{o})=\frac{q}{q-1}$. Rewrite the second integral as
$$\int_\mathfrak{m} \psi(\frac{y}{\varpi^{N}})\,\int_{\mathfrak{o}^\times}\chi(x)\,\psi(\frac{-xy}{\varpi^{N}})\,dx\,dy$$
This is $0$ since
$$\int_{\mathfrak{o}^\times}\chi(x)\,\psi(\frac{-xy}{\varpi^{N}})\,dx = \sum_{a \in \frac{\mathfrak{o^\times}}{1+\varpi^{N-1}\mathfrak{o}}} \chi(a) \, \int_{1+\varpi^{N-1} \mathfrak{o}} \psi(\frac{-xy}{\varpi^{N}})\,dx = 0$$
because $\psi(\frac{x}{\varpi^{N}})$ is trivial on $\varpi^{N}\mathfrak{o}$ and non-trivial on $\varpi^{N-1}\mathfrak{o}$. The integral in $y^\prime$ is:
$$\int_{k^\times} \Bigl(\int_{\mathfrak{o}^\times}{\psi}(xy^\prime t)\,\overline{\chi}(y^\prime t)\,dt \Bigr)\, |y^\prime|^{1-s}\,\, \overline{W}\bigl(\begin{smallmatrix}y^\prime & 0\\ 0 & 1 \end{smallmatrix}\bigr)\,dy^\prime,\,\,\,\,\,\,\,\,\, t \in \mathfrak{o}^\times$$
the conjugate of the integral in $y$. Thus, by replacing $y\eta$ with $u$, and $x$ with $m$, the integrals over $\mathfrak{o}^\times$ in $y$ and $y^\prime$ are
$$\Bigl|\int_{\mathfrak{o}^\times}\overline{\psi}(xy\eta)\,\chi(y\eta)\,d\eta \Bigr|^2\;=\; \Bigl|\int_{\mathfrak{o}^\times}\overline{\psi}(\frac{mu}{\varpi^{N}})\,\chi(u)\,du \Bigr|^2\,=\,\frac{q^{2-N}}{(q-1)^2}$$
So the entire nondecoupled local integral becomes
$$\frac{q^{2-N}}{(q-1)^2} \Bigl[\int_{k^\times} |y|^s\,W\bigl(\begin{smallmatrix}y & 0\\ 0 & 1 \end{smallmatrix}\bigr)\,dy \, \cdot \, \int_{k^\times} |y^\prime|^{1-s}\,\, \overline{W}\bigl(\begin{smallmatrix}y^\prime & 0\\ 0 & 1 \end{smallmatrix}\bigr)\,dy^\prime\, \cdot \, \int_k \varphi\bigl(\begin{smallmatrix}1 & x\\ 0 & 1 \end{smallmatrix}\bigr)\, dx \Bigl]$$

\noindent
Recall that the integral is zero unless $\ord(y)=-\ord(x)-N$. So rewrite the integral as:\\

$\dfrac{q^{2-N}}{(q-1)^2}\, \cdot \,\displaystyle\int_k \varphi\bigl(\begin{smallmatrix}1 & x\\ 0 & 1 \end{smallmatrix}\bigr)\, \cdot \, \displaystyle\int_{\ord(y)=-\ord(x)-N} |y|^s\,W\bigl(\begin{smallmatrix}y & 0\\ 0 & 1 \end{smallmatrix}\bigr)\,dy \, \cdot \, \displaystyle\int_{\ord(y^\prime)=-\ord(x)-N} |y^\prime|^{1-s}\,\overline{W}\bigl(\begin{smallmatrix}y^\prime & 0\\ 0 & 1 \end{smallmatrix}\bigr)\,dy^\prime\, \, dx$\\

$=\dfrac{q^{2-N}}{(q-1)^2}\, \cdot \,\displaystyle\int_k \varphi\bigl(\begin{smallmatrix}1 & x\\ 0 & 1 \end{smallmatrix}\bigr)\, \cdot \, \displaystyle\int_{\ord(y)=-\ord(x)-N} |y|\, |W\bigl(\begin{smallmatrix}y & 0\\ 0 & 1 \end{smallmatrix}\bigr) |^2\,dy\, dx$\\

\noindent
Since $\ord(y)=-\ord(x)-N$, then $y$ can be written as
$$y=\frac{t}{\varpi^{N} x}, \,\,\, t \in \mathfrak{o}^\times$$
So the entire integral is:
$$\frac{q^{2-N}}{(q-1)^2}\, \cdot \,\int_k \varphi\bigl(\begin{smallmatrix}1 & x\\ 0 & 1 \end{smallmatrix}\bigr)\, \cdot \, |\frac{1}{\varpi^{N} x} |\, |W\bigl(\begin{smallmatrix} \frac{1}{\varpi^{N} x} & 0\\ 0 & 1 \end{smallmatrix}\bigr) |^2\,dx$$
Now $y \to W\bigl(\begin{smallmatrix} y & 0\\ 0 & 1 \end{smallmatrix}\bigr)$ is supported on $\mathfrak{o}\cap k^\times$, so $\ord(x) \leq -N$. Then $x \not\in \mathfrak{o}$. Thus, we integrate over $k^\times$, and by a change in Haar measure, the integral becomes\\

$\dfrac{q^{2-N}}{(q-1)^2}\, \cdot \,\displaystyle\int_k |x|^{-w^\prime}\, \cdot \, |\dfrac{1}{\varpi^{N} x} |\, |W\bigl(\begin{smallmatrix} \frac{1}{\varpi^{N} x} & 0\\ 0 & 1 \end{smallmatrix}\bigr) |^2\,dx$\\

$= \dfrac{q^{1-N}}{q-1}\, \cdot \,\displaystyle\int_{k^\times} |x|^{1-w^\prime}\, \cdot \, |\dfrac{1}{\varpi^{N} x} |\, |W\bigl(\begin{smallmatrix} \frac{1}{\varpi^{N} x} & 0\\ 0 & 1 \end{smallmatrix}\bigr) |^2\,dx$\\

\noindent
Invert $x$ to get
$$\frac{q^{1-N}}{q-1}\, \cdot \,\int_{k^\times} |x|^{w^\prime-1}\, \cdot \, |\frac{x}{\varpi^{N}} |\, |W\bigl(\begin{smallmatrix} \frac{x}{\varpi^{N}} & 0\\ 0 & 1 \end{smallmatrix}\bigr) |^2\,dx$$
Replace $x$ by $\varpi^{N}x$ and let $\ord(x)=\ell$ to get\\

$\mathcal{K}_{v_1}(w^\prime,\chi_{v_1})\,=\,\dfrac{q^{1-N}}{q-1}\, \cdot \, q^{-Nw^\prime}\,\cdot \,q^N \, \cdot \,\displaystyle\int_{k^\times} |x|^{w^\prime-1}\, \cdot \, |x| \,\cdot \, |W\bigl(\begin{smallmatrix} x & 0\\ 0 & 1 \end{smallmatrix}\bigr) |^2\,dx$\\

$= \dfrac{q}{q-1}\, \cdot \, q^{-Nw^\prime}\, \cdot \, \displaystyle\sum_{\ell=0}^\infty q^{-\ell w^\prime}\, \cdot \, \dfrac{\alpha^{\ell+1}-\beta^{\ell+1}}{\alpha-\beta}\, \cdot \,\dfrac{\overline{\alpha}^{\ell+1}-\overline{\beta}^{\ell+1}}{\overline{\alpha}-\overline{\beta}}$\\

$=\dfrac{q^{1-Nw^\prime}}{q-1}\, \cdot \, \dfrac{1-|\alpha|^2 |\beta|^2 q^{-2w^\prime}}{(1-|\alpha|^2q^{-w^\prime})(1-|\beta|^2q^{-w^\prime})(1-\overline{\alpha}\beta q^{-w^\prime})(1-\alpha \overline{\beta}q^{-w^\prime})}$\\

\vspace{4mm}
\newpage
\noindent
{\bf 3.1.4 Asymptotics}\\

\noindent
Now
$$\dfrac{q}{q-1}\, \cdot \, \dfrac{1-|\alpha|^2 |\beta|^2 q^{-2w^\prime}}{(1-|\alpha|^2q^{-w^\prime})(1-|\beta|^2q^{-w^\prime})(1-\overline{\alpha}\beta q^{-w^\prime})(1-\alpha \overline{\beta}q^{-w^\prime})}$$
is independent of the conductor $q^N$ of $\chi$, so
$$\mathcal{K}_{v_1}(w^\prime,\chi_{v_1}) \ll (q^N)^{-w^\prime}$$

\noindent
The formulas (5.2) - (5.4) in [Diaconu-Garrett 2008] state that the nondecoupled integral 
$$\mathcal{K}_\infty(s,\chi_0,\chi)=\prod_v \mathcal{K}_v(s,\chi_{0,v},\chi_v)=\prod_v \mathcal{K}_v(s,s^\prime,w,\chi_v)$$
has the following asymptotic formula:\\

\noindent
For $v$ complex,\\

$\mathcal{K}_v(s,s^\prime,w,\chi_v)$\\

$=\pi^{-2s^\prime+1}\, A(s^\prime,w,\mu_1,\mu_2)\,\cdot\,(1+\ell_v^2+4(t+t_v)^2)^{-w}\,\cdot \, \Big[1+O\Bigl( \bigl(\sqrt{1+\ell_v^2+4(t+t_v)^2} \bigr)^{-1} \Bigr) \Bigr]$\\

\noindent
where $A(s^\prime,w,\mu_1,\mu_2)$ is the ratio of products of gamma functions\\

$2^{4w-4s^\prime-4} \dfrac{\Gamma(w+s^\prime+i\mu_1+i \overline{\mu}_2)\Gamma(w+s^\prime-i\mu_1+i \overline{\mu}_2)\Gamma(w+s^\prime+i\mu_1-i \overline{\mu}_2)\Gamma(w+s^\prime-i\mu_1-i \overline{\mu}_2)}{\Gamma(2w+2s^\prime)}$\\

\noindent
and $it_v,\,\ell_v$ are the parameters of the local component $\chi_v$ of $\chi$.\\

\noindent
For $v$ real,\\

$\mathcal{K}_v(s,s^\prime,w,\chi_v)\,=\,B(s^\prime,w,\mu_1,\mu_2)\,\cdot\,(1+|t+t_v|)^{-w}\,\cdot$\\ 
$$\Big[1+O\Bigl( \bigl(1+|t+t_v)|)^{-\frac{1}{2}} \Bigr) \Bigr]$$

\noindent
where $B(s^\prime,w,\mu_1,\mu_2)$ is a similar ratio of products of gamma functions.\\

\vspace{4mm}
\newpage
\noindent
\mbox{{\bf 4.1 SPECTRAL DECOMPOSITION OF THE POINCAR\'E SERIES}}

\noindent
{\bf 4.1.1 Prologue}\\

\indent
In this chapter, we spectrally decompose the Poincar\'e series; this is central to the ideas underlying the integral moments of automorphic $L$-functions on $GL_2$ to prove the meromorphic continuation of the Poincar\'e series. The decomposition consists of a leading (non-$L^2$) term, cuspidal part and continuous part.\\

\vspace{4mm}
\noindent
{\bf 4.1.2 The Cuspidal Part}\\

\indent
Let $F$ be a cuspform on $G_{\mathbb{A}}$ generating a spherical representation locally everywhere, and suppose $F$ corresponds to a spherical vector everywhere locally. The $F^{th}$ (cuspidal) component of the spectral decomposition of the Poincar\'e series is $\langle P \acute{e},F \rangle\! \cdot\! F$. So\\

$\langle P \acute{e},F \rangle\,=\,\displaystyle\int_{Z_\mathbb{A} G_k \backslash G_\mathbb{A}} P \acute{e}(g) \cdot \overline{F}(g)\, dg\,=\, \displaystyle\int_{Z_\mathbb{A} G_k \backslash G_\mathbb{A}} \sum_{\gamma \in M_k \backslash G_k} \varphi(\gamma g)\, \overline{F}(g) \, dg$\\

$= \displaystyle\int_{Z_\mathbb{A} M_k \backslash G_\mathbb{A}} \varphi(g)\,\overline{F}(g) \, dg\,=\, \displaystyle\int_{Z_\mathbb{A} M_k \backslash G_\mathbb{A}} \varphi(g)\,\sum_{\xi \in Z_k \backslash M_k} \overline{W}_F(\xi g)\, dg\,=\, \displaystyle\int_{Z_\mathbb{A}\backslash G_\mathbb{A}} \varphi(g)\,\overline{W}_F(g)\,dg$\\

$=\displaystyle\prod_{v<\infty,\, v\neq v_1} \displaystyle\int_{Z_v\backslash G_v} \varphi_v(g)\,\overline{W}_{F,v}(g)\,dg\, \cdot \, \displaystyle\int_{Z_{v_1}\backslash G_{v_1}} \varphi_{v_1}(g)\,\overline{W}_{F,{v_1}}(g)\,dg\, \cdot $\\
$$\prod_{v|\infty} \int_{Z_v\backslash G_v} \varphi_v(g)\,\overline{W}_{F,v}(g)\,dg$$ 

\noindent
Suppress $v$. At finite $v \neq v_1$, by Iwasawa decomposition and right $K_v$-invariance, 
$$\int_{Z \backslash MN} \varphi(mn)\,\overline{W}_F(mn)\,dm\,dn$$
Further, with $Z\backslash MN \approx HN$ and $$W_F(mn)=\psi(mnm^{-1})\,W_F(m)$$ we get
$$\int_H\, \int_N \chi_0(m)\,\varphi(n)\,\overline{\psi}(mnm^{-1})\,\overline{W}_F(m)\,dm\,dn$$ 
Again, for $m$ in the support of $W_F$ and $n \in N \cap K$
$$\int_N \varphi(n)\,\overline{\psi}(mnm^{-1})\,dn\,=\,1$$
So the integral becomes
$$\int_H \chi_0(m)\, \overline{W}_F(m)\,dm \,= \,\int_{k^\times} |y|^{s^\prime}\, \overline{W}\bigl(\begin{smallmatrix}y & 0\\ 0 & 1 \end{smallmatrix}\bigr)\, dy\,=\, (1-\overline{\alpha}q^{-s^\prime})^{-1}\,(1-\overline{\beta}q^{-s^\prime})^{-1}\,= \,L_v(s^\prime+\frac{1}{2},\overline{F})$$

For $v=v_1$, $\langle P \acute{e},F \rangle$ unwinds to
$$\int_H\, \int_N \chi_0(m)\,\varphi(n)\,\overline{\psi}(mnm^{-1})\,\overline{W}_F(m)\,dm\,dn$$ 
Now $$mnm^{-1}=\bigl(\begin{smallmatrix}1 & xy\\ 0 & 1 \end{smallmatrix}\bigr)$$
So the integral becomes
$$\int_k\, \int_{k^\times} \overline{\psi}(xy)\,|y|^{s^\prime}\,\, \overline{W}\bigl(\begin{smallmatrix}y & 0\\ 0 & 1 \end{smallmatrix}\bigr)\,\varphi(x)\,dy\,dx$$
We will first consider the integral in $y$
$$\int_{k^\times} \overline{\psi}(xy)\,|y|^{s^\prime}\,\overline{W}\bigl(\begin{smallmatrix}y & 0\\ 0 & 1 \end{smallmatrix}\bigr)\,dy$$
where $\psi$ is the standard additive character trivial on $\mathfrak{o}$ and non-trivial on $\varpi^{-1}\mathfrak{o}$ for absolutely unramified $v$. The spherical Whittaker function is defined by
\begin{equation*}
W\bigl(\begin{smallmatrix}y & 0\\ 0 & 1 \end{smallmatrix}\bigr)=\\
\begin{cases}
\dfrac{\alpha^{n+1}-\beta^{n+1}}{\alpha-\beta}  & \text{($n \geq 0$)}\\
0  &\text{(otherwise)}
\end{cases}
\end{equation*}

\noindent
where $\overline{W}_F=W_{\overline{F}}=W$, where $\alpha,\beta$ are Satake parameters of $\overline{F}$, and where $\ord(y)=n$. For $x \in \mathfrak{o}$ (i.e for $\ord(x) \geq 0$), since $\psi$ is trivial on $\mathfrak{o}$, the integral is
$$\sum_{n=0}^\infty q^{-ns^\prime}\,\frac{\alpha^{n+1}-\beta^{n+1}}{\alpha-\beta}\,=\,(1-\alpha q^{-s^\prime})^{-1}\,(1-\beta q^{-s^\prime})^{-1}\,=\,L_v(s^\prime+\frac{1}{2},\overline{F})$$
\\
For $x \not\in \mathfrak{o}$ (i.e for $\ord(x) < 0$), we first evaluate
$$\int_{\mathfrak{o^\times}} \overline{\psi}(xy)\, dy$$
Let $\ord(x)=m$. Write
$$y=\varpi^nt,\,\,\,\,\,x=\varpi^m\eta,\,\,\,\,\,\,t,\eta \in \mathfrak{o^\times}$$
Then, replacing $t\eta$ by $u$, the integral becomes
$$\int_{\mathfrak{o^\times}} \overline{\psi}(\varpi^{m+n}u)\, du$$
Now
$$\int_{\mathfrak{o^\times}} \overline{\psi}(\varpi^{\ell}u)\, du=\int_{\mathfrak{o}} \overline{\psi}(\varpi^{\ell}u)\, du-\int_{\mathfrak{m}} \overline{\psi}(\varpi^{\ell}u)\, du$$
For $\ell \geq 0$, the integrand is $1$ so
\newcommand{\meas}{meas}
$$\int_{\mathfrak{o^\times}} \overline{\psi}(\varpi^{\ell}u)\, du=\meas(\mathfrak{o^\times})=1$$
For $\ell=-1$, $\overline{\psi}$ is non-trivial on $\mathfrak{o}$ and trivial on $\mathfrak{m}$, so
$$\int_{\mathfrak{o^\times}} \overline{\psi}(\varpi^{\ell}u)\, du=-\meas(\mathfrak{m})=\,-\,\frac{1}{q-1}$$
For $\ell \leq -2$, $\overline{\psi}$ is non-trivial on $\mathfrak{o}$ and on $\mathfrak{m}$, so
$$\int_{\mathfrak{o^\times}} \overline{\psi}(\varpi^{\ell}u)\, du= 0$$
So keeping in mind that $\ord(y)=n$ and $\ord(x)=m$,
\begin{align*}
\int_{\mathfrak{o^\times}} \overline{\psi}(xy)\,dy=\int_{\mathfrak{o^\times}} \overline{\psi}(\varpi^{m+n}u)\, du=
\begin{cases}
1  &\text{(for $\ord(y) \geq -\ord(x)$)}\\
-\frac{1}{q-1} &\text{(for $\ord(y)=-\ord(x)-1$)}\\
0 &\text{(otherwise)}
\end{cases}
\end{align*}
So, for $x \not\in \mathfrak{o}$,
$$\int_{k^\times} \overline{\psi}(xy)\,|y|^{s^\prime}\,W\bigl(\begin{smallmatrix}y & 0\\ 0 & 1 \end{smallmatrix}\bigr)\,dy$$
$$=\int_{\ord(y) \geq -\ord(x)} |y|^{s^\prime}\,W\bigl(\begin{smallmatrix}y & 0\\ 0 & 1 \end{smallmatrix}\bigr)\,dy\,-\,\frac{1}{q-1} \int_{\ord(y)=-\ord(x)-1} |y|^{s^\prime}\,W\bigl(\begin{smallmatrix}y & 0\\ 0 & 1 \end{smallmatrix}\bigr)\,dy$$
\\
\noindent
Now the whole integral is
$$\int_{k}\, \int_{k^\times} \overline{\psi}(xy)\,|y|^{s^\prime}\,W\bigl(\begin{smallmatrix}y & 0\\ 0 & 1 \end{smallmatrix}\bigr)\, \varphi(x)\,dy\, dx$$
Again, the sub-integral over $x \in \mathfrak{o}$ evaluates to
$$(1-\alpha q^{-s^\prime})^{-1}\,(1-\beta q^{-s^\prime})^{-1}=L_v(s^\prime+\frac{1}{2},f)$$
The sub-integral over $x \not \in \mathfrak{o}$ becomes\\

$\displaystyle\int_{\ord(x)<0} |x|^{-w^\prime}\,\displaystyle\int_{\ord(y) \geq -\ord(x)} |y|^{s^\prime}\,W\bigl(\begin{smallmatrix}y & 0\\ 0 & 1 \end{smallmatrix}\bigr)\,dy\,dx\,$\\
$$-\,\frac{1}{q-1}\int_{\ord(x)<0} |x|^{-w^\prime}\, \int_{\ord(y)=-\ord(x)-1} |y|^{s^\prime}\,W\bigl(\begin{smallmatrix}y & 0\\ 0 & 1 \end{smallmatrix}\bigr)\,dy\,dx$$

\noindent
First,\\ 

$\displaystyle\int_{\ord(x)<0} |x|^{-w^\prime}\,\displaystyle\int_{\ord(y) \geq -\ord(x)} |y|^{s^\prime}\,W\bigl(\begin{smallmatrix}y & 0\\ 0 & 1 \end{smallmatrix}\bigr)\,dy\,dx$\\

$=\dfrac{q-1}{q}\displaystyle\int_{\ord(x)<0} |x|^{1-w^\prime}\,\displaystyle\int_{\ord(y) \geq -\ord(x)} |y|^{s^\prime}\,W\bigl(\begin{smallmatrix}y & 0\\ 0 & 1 \end{smallmatrix}\bigr)\,dy\,dx$\\

$= \dfrac{q-1}{q} \displaystyle\sum_{m=1}^\infty q^{m(1-w^\prime)}\,\cdot\, \displaystyle\sum_{n \geq -m} q^{-ns^\prime}\dfrac{\alpha^{n+1}-\beta^{n+1}}{\alpha-\beta}$\\

$= \dfrac{q-1}{q} \displaystyle\sum_{m=1}^\infty q^{m(1-w^\prime)}\,\cdot\,(\alpha-\beta)^{-1}\bigl[\dfrac{\alpha^{1-m}q^{s^\prime m}}{1-\alpha q^{-s^\prime}}- \dfrac{\beta^{1-m}q^{s^\prime m}}{1-\beta q^{-s^\prime}}\bigr]$\\

$= \dfrac{q-1}{q}\,\cdot\,(\alpha-\beta)^{-1}\bigl[(1-\alpha q^{-s^\prime})^{-1}\,\alpha\,\displaystyle\sum_{m=1}^\infty (\alpha^{-1}q^{1-w^\prime+s^\prime})^m\, - $\\
$$(1-\beta q^{-s^\prime})^{-1}\,\beta\,\sum_{m=1}^\infty (\beta^{-1}q^{1-w^\prime+s^\prime})^m\bigr]$$

$= \dfrac{q-1}{q}\,\cdot\,(\alpha-\beta)^{-1}\bigl[\dfrac{q^{1-w^\prime+s^\prime}}{(1-\alpha q^{-s^\prime})(1-\alpha^{-1}q^{1-w^\prime+s^\prime})}\, - \, \dfrac{q^{1-w^\prime+s^\prime}}{(1-\beta q^{-s^\prime})(1-\beta^{-1}q^{1-w^\prime+s^\prime})}\bigr]$\\

$=\dfrac{(q-1)(q^{-w^\prime}-\dfrac{q^{1-2w^\prime+2s^\prime}}{\alpha\beta})}{(1-\alpha q^{-s^\prime})\,(1-\beta q^{-s^\prime})\,(1-\alpha^{-1}q^{1-w^\prime+s^\prime})\,(1-\beta^{-1}q^{1-w^\prime+s^\prime})}$\\

$=L_{v_1}(s^\prime+\frac{1}{2},f)\,\cdot\,\dfrac{(q-1)(q^{-w^\prime}-\dfrac{q^{1-2w^\prime+2s^\prime}}{\alpha\beta})}{(1-\alpha^{-1}q^{1-w^\prime+s^\prime})\,(1-\beta^{-1}q^{1-w^\prime+s^\prime})}$\\

For $\ord(y)=-\ord(x)-1$, write  $y$ as
$$y=\frac{t}{\varpi x}\,\,\,\, (t \in \mathfrak{o}^\times)$$
Then
$$\frac{1}{q-1}\, \cdot \,\int_{\ord(x)<0} |x|^{-w^\prime}\, \cdot \, \int_{\ord(y)=-\ord(x)-1} |y|^{s^\prime}\, W\bigl(\begin{smallmatrix}y & 0\\ 0 & 1 \end{smallmatrix}\bigr) \,dy\, dx$$
$$=\frac{1}{q-1}\, \cdot \,\int_{\ord(x)<0} |x|^{-w^\prime}\, \cdot \, |\frac{1}{\varpi x} |^{s^\prime}\, W\bigl(\begin{smallmatrix} \frac{1}{\varpi x} & 0\\ 0 & 1 \end{smallmatrix}\bigr) \,dx$$
Now $y \to W\bigl(\begin{smallmatrix} y & 0\\ 0 & 1 \end{smallmatrix}\bigr)$ is supported on $\mathfrak{o}\cap k^\times$. Thus, integrate over $k^\times$, and by changing to multiplicative Haar measure, the integral becomes\\
$$\frac{q-1}{q}\,\cdot\,\frac{1}{q-1}\, \cdot \,\int_{\ord(x)<0} |x|^{1-w^\prime}\, \cdot \, |\frac{1}{\varpi x} |^{s^\prime}\, W\bigl(\begin{smallmatrix} \frac{1}{\varpi x} & 0\\ 0 & 1 \end{smallmatrix}\bigr) \,dx$$
Invert $x$ to obtain
$$\frac{1}{q}\, \cdot \,\int_{\ord(x)>0}|x|^{w^\prime-1}\, \cdot \, |\frac{x}{\varpi}|^{s^\prime}\, W\bigl(\begin{smallmatrix} \frac{x}{\varpi} & 0\\ 0 & 1 \end{smallmatrix}\bigr) \,dx$$
Replace $x$ by $\varpi x$ and with $\ord(x)=m$\\

$\dfrac{1}{q}\, \cdot \, q^{1-w^\prime}\,\displaystyle\int_{\ord(x) \geq 0} |x|^{w^\prime-1}\,\cdot \, |x|^{s^\prime}\,\cdot \, W\bigl(\begin{smallmatrix} x & 0\\ 0 & 1 \end{smallmatrix}\bigr) \,dx$\\

$=q^{-w^\prime}\, \cdot \, \displaystyle\sum_{m=0}^\infty q^{-m (w^\prime-1+s^{\prime})}\, \cdot \, \frac{\alpha^{m+1}-\beta^{m+1}}{\alpha-\beta}$\\

$=q^{-w^\prime}\,\cdot\, (1-\alpha q^{1-w^\prime-s^\prime})^{-1}\,(1-\beta q^{1-w^\prime-s^\prime})^{-1}\,=\,q^{-w^\prime}\,L(s^\prime+w^\prime-\frac{1}{2},F)$\\

\noindent
So, for $v=v_1$, the $v^{th}$ local factor of $\langle P \acute{e},F \rangle$ is\\
$$\frac{1}{(1-\alpha q^{-s^\prime})\,(1-\beta q^{-s^\prime})}\,+\,\frac{(q-1)(q^{-w^\prime}-\dfrac{q^{1-2w^\prime+2s^\prime}}{\alpha\beta})}{(1-\alpha q^{-s^\prime})\,(1-\beta q^{-s^\prime})\,(1-\alpha^{-1}q^{1-w^\prime+s^\prime})\,(1-\beta^{-1}q^{1-w^\prime+s^\prime})}$$
$\,-\,\dfrac{1}{q^{w^\prime}\,(1-\alpha q^{1-w^\prime-s^\prime})\,(1-\beta q^{1-w^\prime-s^\prime})}$\\

\mbox{$=L_v(s^\prime+\frac{1}{2}, \overline{F})+\dfrac{(q-1)(q^{-w^\prime}-\dfrac{q^{1-2w^\prime+2s^\prime}}{\alpha\beta})}{(1-\alpha^{-1}q^{1-w^\prime+s^\prime})\,(1-\beta^{-1}q^{1-w^\prime+s^\prime})}\,\cdot\,L_v(s^\prime+\frac{1}{2}, \overline{F})\,-\,q^{-w^\prime}L_v(s^\prime+w^\prime-\frac{1}{2}, \overline{F})$}
\\

For infinite $v$, by formulas (4.2) and (4.3) in [Diaconu-Garrett 2008] the $v^{th}$ local factor of $\langle P \acute{e},F \rangle$ is $ \mathcal{G}(\frac{1}{2}+i \overline{\mu}_{F,v};s^\prime,w)$, where up to a constant, \mbox{for $v \approx \mathbb{R}$,}
$$\mathcal{G}_v(s;s^\prime,w)=\pi^{-s^\prime}\,\dfrac{\Gamma(\dfrac{s^\prime+1-s}{2})\,\Gamma(\dfrac{s^\prime+w-s}{2})\,\Gamma(\dfrac{s^\prime+s}{2})\,\Gamma(\dfrac{s^\prime+w+s-1}{2})}{\Gamma(\dfrac{w}{2})\,\Gamma(s^\prime+\dfrac{w}{2})}$$
and at $v \approx \mathbb{C}$,
$$\mathcal{G}_v(s;s^\prime,w)=2\pi^{-2s^\prime}\,\frac{\Gamma(s^\prime+1-s)\,\Gamma(s^\prime+w-s)\,\Gamma(s^\prime+s)\,\Gamma(s^\prime+w+s-1)}{\Gamma(w)\,\Gamma(2s^\prime+w)}$$
Group the archimedean factors as
$$\mathcal{G}_{F_\infty}(s^\prime,w)=\prod_{v|\infty}\mathcal{G}_v(\frac{1}{2}+i \overline{\mu}_{F,s^\prime};s^\prime,w)$$
and let all ambiguous constants be absorbed into $\overline{\rho}_F$. Then, for cuspforms $F$, the cuspidal part of the spectral decomposition of the Poincar\'e series is\\

$\sum_F \langle P \acute{e},F \rangle\! \cdot\! F\,=\,\sum_F \overline{\rho}_F \mathcal{G}_{F_\infty}(s^\prime,w)\,\cdot\, [L_v(s^\prime+\frac{1}{2},\overline{F})+L_{v_1}(s^\prime+\frac{1}{2}, \overline{F})\,+$\\

\mbox{$\dfrac{(q-1)(q^{-w^\prime}-\dfrac{q^{1-2w^\prime+2s^\prime}}{\alpha\beta})}{(1-\alpha^{-1}q^{1-w^\prime+s^\prime})\,(1-\beta^{-1}q^{1-w^\prime+s^\prime})}\,\cdot\,L_{v_1}(s^\prime+\frac{1}{2}, \overline{F})\,-\,q^{-w^\prime}L_{v_1}(s^\prime+w^\prime-\frac{1}{2},\overline{F})]\,\cdot\,F$}\\

There is no residual spectrum since residual automorphic forms on $GL(2)$ are associated to one-dimensional representations which have no Whittaker models.\\
\\

\vspace{4mm}
\noindent
{\bf 4.1.3 The continuous part}\\

\indent
Subtract an Eisenstein series from the Poincar\'e series and denote the resulting function by $P\acute{e}^*$. This function is $L^2$ and has sufficient decay so that it can be integrated against an Eisenstein series (see section 4 in [Diaconu-Garrett 2008]). The leading term is:

$\displaystyle\int_{N_\mathbb{A}} \varphi\,=\,\Bigl(\displaystyle\int_{N_\infty} \varphi_\infty \, \cdot \, \Bigl[\displaystyle\int_{N_{v \neq v_1}} \varphi_{v \neq v_1} \, \cdot \, \displaystyle\int_{N_{v_1}} \varphi_{v_1}\Bigr] \Bigr) \,\cdot\, E_{s^\prime+1,1}$\\
\noindent
where, as in (4.16) in [Diaconu-Garrett 2008] an elementary computation shows
\begin{align*}
\int_{N_v} \varphi_v =
\begin{cases}
\sqrt{\pi} \frac{\Gamma(\frac{w-1}{2})}{\Gamma(\frac{w}{2})} & \text {($v \approx \mathbb{R}$)}\\
2\pi(w-1)^{-1} & \text{($v \approx \mathbb{C}$)}
\end{cases}
\end{align*}
Now for $v \neq v_1$,
$$\int_{N_v} \varphi_v \, dn\,=\, \int_{k_v} 1\, dx\,=\,1$$
and
$\displaystyle\int_{N_{v_1}}\varphi_{v_1}\,dn\,=\,\displaystyle\int_{x \in \mathfrak{o}_v} 1\,dx\,+\, \displaystyle\int_{x \not \in \mathfrak{o}_v} |x|^{-w^\prime}\,dx$\\

$=1\,+\,\dfrac{q-1}{q}\,\cdot\, \displaystyle\sum_{m=1}^\infty (q^m)^{1-w^\prime}\,=\,1 \,+\, \dfrac{q-1}{q}\,\cdot\,\dfrac{q^{1-w^\prime}}{1-q^{1-w^\prime}}\,=\,\dfrac{1-q^{-w^\prime}}{1-q^{1-w^\prime}}$\\

So the leading term is
$$\int_{N_\infty} \varphi_\infty \,\cdot\, \dfrac{1-q^{-w^\prime}}{1-q^{1-w^\prime}}\,\cdot\,E_{s^\prime+1,1}$$
The continuous part of the spectral decomposition of $P\acute{e}$ is
$$\frac{1}{4\pi i \kappa} \sum_\chi \int_{Re(s)=\frac{1}{2}} \langle P\acute{e}^*,E_{s,\chi} \rangle \cdot E_{s,\chi} \, ds\,\,\,\,\,(\mbox{where}\,\,\,\kappa=\meas(\mathbb{J}^1/k^\times)) $$
So the spectral decomposition of the Poincar\'e series is\\

$P\acute{e}=\Bigl(\displaystyle\int_{N_\infty} \varphi_\infty \Bigr)\,\cdot\,\dfrac{1-q^{-w^\prime}}{1-q^{1-w^\prime}} \,\cdot\, E_{s^\prime+1,1}$\\

$+\displaystyle\sum_F \overline{\rho}_F \mathcal{G}_{F_\infty}(s^\prime,w)\,\cdot\,[L_v(s^\prime+\frac{1}{2},\overline{F})+L_{v_1}(s^\prime+\frac{1}{2}, \overline{F})\,+$\\

\mbox{$\dfrac{(q-1)(q^{-w^\prime}-\dfrac{q^{1-2w^\prime+2s^\prime}}{\alpha\beta})}{(1-\alpha^{-1}q^{1-w^\prime+s^\prime})\,(1-\beta^{-1}q^{1-w^\prime+s^\prime})}\,\cdot\,L_{v_1}(s^\prime+\frac{1}{2}, \overline{F})\,-\,q^{-w^\prime}L_{v_1}(s^\prime+w^\prime-\frac{1}{2},\overline{F})]\,\cdot\,F$}\\

$\,+\,\frac{1}{4\pi i \kappa} \displaystyle\sum_\chi \int_{Re(s)=\frac{1}{2}} \langle P\acute{e}^*,E_{s,\chi} \rangle \cdot E_{s,\chi} \, ds $\\

\noindent
As in section 4 in [Diaconu-Garrett 2008],
$$\langle P\acute{e}^*,E_{s,\chi} \rangle=\Bigl( \int_{Z_\infty \backslash G_\infty} \varphi_\infty \,\cdot\, \overline{W}^E_{s,\chi,\infty} \Bigr)\,\cdot\, \Bigl( \prod_{v<\infty}  \int_{Z_v \backslash G_v} \varphi_v(g_v) \,\cdot\, \overline{W}^E_{s,\chi,v}(g_v) \Bigr)\, dg_v$$
where
\begin{align*}
\int_{Z_v \backslash G_v} \varphi_\infty \,\cdot\, \overline{W}^E_{s,\chi,v}=
\begin{cases}
\displaystyle\frac{\mathcal{G}_v(s,s^\prime,w)}{\pi^{-s}\Gamma(s)}  & \text{($v \approx \mathbb{R}$)}\\
\displaystyle\frac{\mathcal{G}_v(s,s^\prime,w)}{2 \pi^{-2s-1}\Gamma(2s)}  & \text{($v \approx \mathbb{C}$)}
\end{cases}
\end{align*}
and for finite $v \neq v_1$, \\

$\displaystyle\int_{Z_v \backslash G_v} \varphi_v(g_v) \,\cdot\, \overline{W}^E_{s,\chi,v}(g_v)$\\
$$=|\mathfrak{d}_v|_v^{\frac{1}{2}}\, \cdot\, \frac{L_v(s^\prime+\overline{s},\overline{\chi}_v)\,\cdot\,L_v(s^\prime+1-\overline{s},\chi_v)}{L_v(2\overline{s},\overline{\chi}_v^2)}\,\cdot\, |\mathfrak{d}_v|_v^{-(s^\prime+1-\overline{s})}\,\cdot\, \overline{\chi}_v(\mathfrak{d}_v)$$
where $\mathfrak{d}$ is the idele with $v^{th}$ component $\mathfrak{d}_v$ at finite place $v$ and component $1$ at archimedean places.\\

\noindent
For finite $v=v_1$,
$$\int_{Z_v \backslash G_v} \varphi_v(g_v) \,\cdot\, \overline{W}^E_{s,\chi,v}(g_v)=\int_k\,\int_{k^\times}  |y|^{s^\prime}\,\overline{\psi}(xy)\,\overline{W}_{s,\chi}^E \bigl(\begin{smallmatrix}y & 0\\ 0 & 1 \end{smallmatrix}\bigr)\,\cdot\, \varphi(x)\,dy\,dx$$
Define an Eisenstein series by
$$E(g)=\sum_{\lambda \in P_k \backslash G_k} \eta(\lambda g)$$
for $\eta$ left $P_k$-invariant, left $M_k$-invariant and left $N_\mathbb{A}$-invariant. Present the vectors $\eta_v$ in a different form, namely
$$\eta_v (pk)= \left|\frac{a}{d}\right|_v^{s}\,\cdot \,\chi_v\bigl(\frac{a}{d} \bigr)\,\,\,\,\, (\mbox{for}\,\,\, p=(\begin{smallmatrix} a & *\\0 & d \end{smallmatrix}) \in P_v,\,k \in K_v)$$
Let $\phi_v$ be any Schwartz function on $k_v^2$, invariant under $k_v$ and put
$$\eta_v^{\prime}(g)=\chi_v(\det g) |\det g|^s_v \,\cdot\, \int_{k_v^{\times}} \chi_v^2(t) |t|_v^{2s}\,\cdot\, \phi_v(t \cdot e_2 \cdot g)\,dt$$
where $e_2=e_{2,v}$ is the second basis element in $k_v^2$. $\eta_v^{\prime}$ has the same left $P_v$-equivariance as $\eta_v$:
$$\eta_v^{\prime} \bigl((\begin{smallmatrix} a & *\\0 & d \end{smallmatrix})\,\cdot\,g \bigl)= \left|\frac{a}{d}\right|_v^{s}\,\cdot \,\chi_v\bigl(\frac{a}{d} \bigr)\,\cdot\,\eta_v^\prime(g)$$
For $\phi_v$ invariant under $K_v$, the function $\eta_v^\prime$ is right $K_v$-invariant. So as in Appendix 2 in [Diaconu-Garrett 2008],\\
$$\eta_v^\prime(g)=\eta_v^\prime(1)\,\cdot\,\eta_v(g)\,\,\,\,(\mbox{since}\,\,\,\eta_v(1)=1)$$
and 
$$\eta_v^\prime(1)=\int_{k_v^2} \chi_v^2(t)\, |t|_v^{2s}\,\cdot\, \phi_v(t \cdot e_2 \cdot 1)\,dt= \zeta_v(2s,\chi^2,\phi(0,*))$$

\noindent
Thus, it suffices to compute the local Mellin transform of\\

$\eta_v^\prime(1)\,\cdot\,W_{s,\chi,v}^E(m)\,= \,\displaystyle\int_{N_v} \overline{\psi}(n)\,\cdot\,\eta_v^\prime(w_0 nm)\,dn$\\

$=\chi(y)|y|^s\,\cdot\,\displaystyle\int_{N_v} \overline{\psi}(n)\,\displaystyle\int_{k_v^{\times}} \chi_v^2(t)\, |t|_v^{2s}\,\cdot\, \phi_v(t \cdot e_2 \cdot w_0 \cdot nm)\,dt\,dn$\\

$=\chi(y)|y|^s\,\cdot\,\displaystyle\int_{k_v} \overline{\psi}(x^\prime)\,\displaystyle\int_{k_v^{\times}} \chi_v^2(t)\, |t|_v^{2s}\,\cdot\, \phi_v(tx^\prime,ty)\,dt\,dx^\prime\,\,\,\,\,(\mbox{with}\,\,\,m=(\begin{smallmatrix} y & 0\\0 & 1 \end{smallmatrix}))$\\

\noindent
At finite primes, take
\newcommand{\ch}{ch} 
$$\phi(t,x^\prime)=\ch_{\mathfrak{o}_v}(t)\,\cdot\,\ch_{\mathfrak{o}_v}(x^\prime)\,\,\,\,\,\,\,\,\,(\ch_X=\,\mbox{characteristic function of a set}\,\,X)$$
Then
$$\eta_v^\prime(1)=\zeta_v(2s,\chi^2,\ch_{\mathfrak{o}_v})=L_v(2s,\chi^2)$$
and\\

$\eta_v^\prime(1)\,\cdot\,W_{s,\chi,v}^E (\begin{smallmatrix} y & 0\\0 & 1 \end{smallmatrix})\,=\,\chi(y)|y|^s\,\cdot\,\displaystyle\int_{k_v} \overline{\psi}(x^\prime)\,\ch_{\mathfrak{o}_v}(tx^\prime)\,\cdot \,\displaystyle\int_{k_v^{\times}} \chi_v^2(t)\, |t|_v^{2s}\,\ch_{\mathfrak{o}_v}(ty)\,dt\,dx^\prime$\\

$=\chi(y)|y|^s\,\meas(\mathfrak{o}_v)\,\displaystyle\int_{k_v^{\times}} \ch_{\mathfrak{o}_v^*}(\frac{1}{t})\,\chi_v^2(t)\, |t|_v^{2s-1}\,\ch_{\mathfrak{o}_v}(ty)\,dt$\\

$=|\mathfrak{d}_v|^{\frac{1}{2}}\,\cdot\, \chi(y)|y|^s\,\displaystyle\int_{k_v^{\times}} \ch_{\mathfrak{o}_v^*}(\frac{1}{t})\,\chi_v^2(t)\, |t|_v^{2s-1}\,\ch_{\mathfrak{o}_v}(ty)\,dt$\\

\noindent
where $\mathfrak{d}_v \in k_v^\times$ is such that $(\mathfrak{o}_v^*)^{-1}=\mathfrak{d}_v\,\cdot\,\mathfrak{o}_v$.\\

\noindent
So, omitting $|\mathfrak{d}_v|^{\frac{1}{2}}$ for now,\\

$\displaystyle\int_k\,\displaystyle\int_{k^\times}  |y|^{s^\prime}\,\overline{\psi}(xy)\,\overline{W}_{s,\chi}^E \bigl(\begin{smallmatrix}y & 0\\ 0 & 1 \end{smallmatrix}\bigr)\,\cdot\, \varphi(x)\,dy\,dx$\\

$=\displaystyle\int_k\,\displaystyle\int_{k^\times}  |y|^{s^\prime}\,\overline{\psi}(xy)\, \cdot \, \bigl(\chi(y)|y|^s\,\displaystyle\int_{k_v^{\times}} \ch_{\mathfrak{o}_v^*}(\frac{1}{t})\,\chi_v^2(t)\, |t|_v^{2s-1}\,\ch_{\mathfrak{o}_v}(ty)\,dt \bigr)\, \cdot \,\varphi(x)\,dy\,dx$\\

\noindent
Consider the integrals in $y$ and $t$. Replace $y$ by $\frac{y}{t}$ to get\\

$\displaystyle\int_{k^\times}\, \displaystyle\int_{k^\times} \overline{\psi}(x \frac{y}{t})\, \cdot \, \chi(y)\,|y|^{s+s^\prime}\,\ch_{\mathfrak{o}_v}(y)\, \cdot \ch_{\mathfrak{o}_v^*}(\frac{1}{t})\,\cdot\,\chi(t)\,|t|^{s-1-s^\prime}\,dt \, dy$\\

\noindent
Replace $t$ by $\frac{1}{t}$ to get\\

$\displaystyle\int_{k^\times}\, \displaystyle\int_{k^\times} \overline{\psi}(xyt)\, \cdot \, \chi(y)\,|y|^{s+s^\prime}\,\ch_{\mathfrak{o}_v}(y)\, \cdot \ch_{\partial_v \mathfrak{o}_v}(t)\,\cdot\,\overline{\chi}(t)\,|t|^{s^\prime+1-s}\,dt \, dy$\\

\noindent
First consider the integral in $y$:
$$\int_{k^\times} \overline{\psi}(xty)\, \cdot \, \chi(y)\,|y|^{s+s^\prime}\,\ch_{\mathfrak{o}_v}(y)\,dy $$
\noindent
For $x \in \mathfrak{o}^\times$, $\psi$ is trivial on $\mathfrak{o}$, so we get
$$\int_{\mathfrak{o}^\times} \chi(y)\,|y|^{s+s^\prime}\,dy\,\cdot\,\int_{k^\times}\ch_{\partial_v \mathfrak{o}_v}(t)\,\cdot\,\overline{\chi}(t)\,|t|^{s^\prime+1-s}\,dt$$
$$=L_{v_1}(s+s^\prime,\chi)\,\cdot\,L_{v_1}(s^\prime+1-s,\overline{\chi})\,\cdot\,|\mathfrak{d}_v|^{-(s^\prime+1-s)}\,\chi(\mathfrak{d}_v)$$

\noindent
For $x \not \in \mathfrak{o}^\times$,
\begin{align*}
\int_{\mathfrak{o^\times}} \overline{\psi}(xty)\,dy=
\begin{cases}
1  &\text{($\ord(y) \geq -\ord(x)-\ord(t)$)}\\
-\frac{1}{q-1} &\text{$(\ord(y)=-\ord(x)-\ord(t)-1$)}\\
0 &\text{(otherwise)}
\end{cases}
\end{align*}
So\\

$\displaystyle\int_{k^\times} \overline{\psi}(xty)\, \cdot \, \chi(y)\,|y|^{s+s^\prime}\,\ch_{\mathfrak{o}_v}(y)\,dy \,=\,\displaystyle\int_{\mathfrak{o}^\times} \overline{\psi}(xty)\, \cdot \, \chi(y)\,|y|^{s+s^\prime}\,dy $\\

$=\displaystyle\int_{\ord(y) \geq -\ord(x)-\ord(t)} \chi(y)\,|y|^{s+s^\prime}\,dy\,-\,\frac{1}{q-1}\displaystyle\int_{\ord(y) = -\ord(x)-\ord(t)-1} \chi(y)\,|y|^{s+s^\prime}\,dy$\\

\noindent
The entire integral in $t$ and $y$ is:\\

$\displaystyle\int_{x \not \in \mathfrak{o}} \varphi(x)\, \displaystyle\int_{k^\times}\ch_{\mathfrak{d}_v \mathfrak{o}_v}(t)\,\cdot\,\overline{\chi}(t)\,|t|^{s^\prime+1-s}\,dt\, \cdot \,$\\

$\Bigl[\displaystyle\int_{\ord(y) \geq -\ord(x)-\ord(t)} \chi(y)\,|y|^{s+s^\prime}\,dy\,-\,\frac{1}{q-1}\displaystyle\int_{\ord(y) = -\ord(x)-\ord(t)-1} \chi(y)\,|y|^{s+s^\prime}\,dy\Bigr]$\\

\noindent
First take\\

$\displaystyle\int_{x \not \in \mathfrak{o}} \varphi(x)\, \displaystyle\int_{k^\times}\ch_{\mathfrak{d}_v \mathfrak{o}_v}(t)\,\cdot\,\overline{\chi}(t)\,|t|^{s^\prime+1-s}\,dt\, \cdot \,\displaystyle\int_{\ord(y) \geq -\ord(x)-\ord(t)} \chi(y)\,|y|^{s+s^\prime}\,dy\,dx$\\

$=\dfrac{q-1}{q}\,\displaystyle\sum_{m=1}^\infty (q^{m})^{1-w^\prime}\,\cdot\,\displaystyle\sum_{n \geq -m-r}^\infty (q^{-n})^{s+s^\prime}\,\cdot\,\displaystyle\int_{k^\times}\ch_{\mathfrak{d}_v \mathfrak{o}_v}(t)\,\cdot\,\overline{\chi}(t)\,|t|^{s^\prime+1-s}\,dt$\\

\noindent
(where $\ord(t)=r$ and $\chi(y)$ is omitted for now)\\

$=\dfrac{q-1}{q}\,\displaystyle\sum_{m=1}^\infty \dfrac{(q^{1-w^\prime})^m\,\cdot\,(q^{s+s^\prime})^m\,\cdot\, (q^r)^{s+s^\prime}}{1-q^{-(s+s^\prime)}}\,\cdot\,\displaystyle\int_{k^\times}\ch_{\mathfrak{d}_v \mathfrak{o}_v}(t)\,\cdot\,\overline{\chi}(t)\,|t|^{s^\prime+1-s}\,dt$\\

$=\dfrac{(q-1)\,q^{1-w^\prime+s+s^\prime}}{q(1-q^{1-w^\prime+s+s^\prime})(1-q^{-s-s^\prime})}\,\cdot\,\displaystyle\int_{k^\times}\ch_{\mathfrak{d}_v \mathfrak{o}_v}(t)\,\cdot\,\overline{\chi}(t)\,|t|^{s^\prime+1-s}\,|t|^{-s-s^\prime}\,dt$\\

$=\dfrac{(q-1)\,q^{-w^\prime+s+s^\prime}}{(1-q^{1-w^\prime+s+s^\prime})(1-q^{-s-s^\prime})}\,\cdot\,L_{v_1}(1-2s,\overline{\chi})$\\

\noindent
Next we take\\

$-\dfrac{1}{q-1}\,\displaystyle\int_{x \not \in \mathfrak{o}} \varphi(x)\, \displaystyle\int_{k^\times}\ch_{\mathfrak{d}_v \mathfrak{o}_v}(t)\,\cdot\,\overline{\chi}(t)\,|t|^{s^\prime+1-s}\,dt\, \cdot \,\displaystyle\int_{\ord(y) = -\ord(x)-\ord(t)-1}  \chi(y)\,|y|^{s+s^\prime}\,dy\,dx$\\

\noindent
Since $\ord(y)=-\ord(x)- \ord(t)-1$, $y$ can be written as
$$y=\frac{1}{\varpi t x}$$
So the entire integral becomes an integral in $t$ and $x$ as follows:\\

$-\dfrac{1}{q-1}\,\displaystyle\int_{x \not \in \mathfrak{o}} \varphi(x)\, \displaystyle\int_{k^\times}\ch_{\mathfrak{d}_v \mathfrak{o}_v}(t)\,\cdot\,\overline{\chi}(t)\,|t|^{s^\prime+1-s}\,dt\, \cdot \,\displaystyle\int_{\ord(y) = -\ord(x)-\ord(t)-1}  \chi(y)\,|y|^{s+s^\prime}\,dy\,dx$\\

$=-\dfrac{1}{q-1}\,\cdot \, \dfrac{q-1}{q}\,\displaystyle\int_{x \not \in \mathfrak{o}} |x|^{1-w^\prime}\, \displaystyle\int_{k^\times}\ch_{\mathfrak{d}_v \mathfrak{o}_v}(t)\,\cdot\,\overline{\chi}(t)\,|t|^{s^\prime+1-s}\,dt\, \cdot \,|\frac{1}{\varpi x t}|^{s+s^\prime}\,dx$\\

$=-q^{s+s^\prime-1}\,\displaystyle\int_{|x|>1} |x|^{1-w^\prime-s-s^\prime}\, \displaystyle\int_{k^\times}\ch_{\mathfrak{d}_v \mathfrak{o}_v}(t)\,\cdot\,\overline{\chi}(t)\,|t|^{1-2s}\,dt\,dx$\\

$=-q^{s+s^\prime-1}\,\displaystyle\sum_{m=1}^\infty (q^{m})^{1-w^\prime-s-s^\prime}\,\cdot\,L_{v_1}(1-2s,\overline{\chi})$\\

$=-\dfrac{q^{s+s^\prime-1}\,\cdot\,q^{1-w^\prime-s-s^\prime}}{1-q^{1-w^\prime-s-s^\prime}}\,\cdot\,L_{v_1}(1-2s,\overline{\chi})$\\

$=-\dfrac{q^{-w^\prime}}{1-q^{1-w^\prime-s-s^\prime}}\,\cdot\,L_{v_1}(1-2s,\overline{\chi})$\\

\noindent
Thus adding up we get\\

$$L_{v_1}(1-2s,\overline{\chi})\,\cdot\,\Bigl[\dfrac{(q-1)\,q^{-w^\prime+s+s^\prime}}{(1-q^{1-w^\prime+s+s^\prime})(1-q^{-s-s^\prime})}\,-\,\dfrac{q^{-w^\prime}}{1-q^{1-w^\prime-s-s^\prime}}\Bigr]$$
\noindent
Thus at finite primes $v=v_1$, the integral evaluates to:\\

$L_{v_1}(s+s^\prime,\chi)\,\cdot\,L_{v_1}(s^\prime+1-s,\overline{\chi})\,\cdot\,|\mathfrak{d}_{v_1}|^{-(s^\prime+1-s)}\,\chi(\mathfrak{d}_{v_1})\,+\,$\\

$L_{v_1}(1-2s,\overline{\chi})\,\cdot\,\Bigl[\dfrac{(q-1)\,q^{-w^\prime+s+s^\prime}}{(1-q^{1-w^\prime+s+s^\prime})(1-q^{-s-s^\prime})}\,-\,\dfrac{q^{-w^\prime}}{1-q^{1-w^\prime-s-s^\prime}}\Bigr]\,\cdot\, |\mathfrak{d}_{v_1}|^{-(1-2s)}\,\chi(\mathfrak{d}_{v_1})$\\

\noindent
Then dividing through by $\eta_v^\prime$ and putting back the measure constant $|\mathfrak{d}_v|^{\frac{1}{2}}$, we get for $v=v_1$,\\

$|\mathfrak{d}_v|^{\frac{1}{2}}\,\cdot\,\displaystyle\int_{k^\times}\,\displaystyle\int_k  |y|^{s^\prime}\,\overline{\psi}(xy)\,\overline{W}_{s,\chi}^E \bigl(\begin{smallmatrix}y & 0\\ 0 & 1 \end{smallmatrix}\bigr)\,\cdot\, \varphi(x)\,dy\,dx$\\

$=\dfrac{L_v(s^\prime+s, \chi)\,\cdot\,L_v(s^\prime+1-s, \overline{\chi})\,\cdot\, |\mathfrak{d}_v|^{-(s^\prime+1-s)}\,\cdot\, |\mathfrak{d}_v|^{\frac{1}{2}}\,\cdot\,\chi(\mathfrak{d}_v)}{L_v(2s,\chi^2)}$\\

$+\, \dfrac{L_v(1-2s,\overline{\chi})\,\cdot\,[\dfrac{(q-1)\,q^{-w^\prime+s+s^\prime}}{(1-q^{1-w^\prime+s+s^\prime})(1-q^{-s-s^\prime})}\,-\,\dfrac{q^{-w^\prime}}{1-q^{1-w^\prime-s-s^\prime}}]\,\cdot\, |\mathfrak{d}_v|^{-(1-2s)}\,\cdot \, |\mathfrak{d}_v|^{\frac{1}{2}}\,\cdot \,\chi(\mathfrak{d}_v)}{L_v(2s,\chi^2)}$\\

\noindent
Replacing $s$ by $1-s$ and $\chi$ by $\overline{\chi}$ we get\\

$|\mathfrak{d}_v|^{\frac{1}{2}}\,\cdot\,\displaystyle\int_{k^\times}\,\displaystyle\int_k  |y|^{s^\prime}\,\overline{\psi}(xy)\,\overline{W}_{1-s,\overline{\chi}}^E \bigl(\begin{smallmatrix}y & 0\\ 0 & 1 \end{smallmatrix}\bigr)\,\cdot\, \varphi(x)\,dy\,dx$\\

$=\dfrac{L_v(s^\prime+s, \chi)\, \cdot \,L_v(s^\prime+1-s,\overline{\chi})\,\cdot\, |\mathfrak{d}_v|^{-(s^\prime+s-\frac{1}{2})}\,\cdot\, \chi(\mathfrak{d}_v)}{L_v(2-2s,\overline{\chi}^2)}$\\

$+\, \dfrac{L_v(2s-1,\overline{\chi})\,\cdot\,[\dfrac{(q-1)\,q^{1-w^\prime+s+s^\prime}}{(1-q^{2-w^\prime-s+s^\prime})(1-q^{-1+s-s^\prime})}\,-\,\dfrac{q^{-w^\prime}}{1-q^{-w^\prime+s-s^\prime}}]\,\cdot\, |\mathfrak{d}_v|^{\frac{3}{2}-2s}\,\cdot \,\chi(\mathfrak{d}_v)}{L_v(2-2s,\chi^2)}$\\

\noindent
So the spectral decomposition of the Poincar\'e series is:\\

$P\acute{e}=\Bigl(\displaystyle\int_{N_\infty} \varphi_\infty \Bigr)\,\cdot\,\dfrac{1-q^{-w^\prime}}{1-q^{1-w^\prime}} \,\cdot\, E_{s^\prime+1,1}$\\

$+\displaystyle\sum_F \overline{\rho}_F \mathcal{G}_{F_\infty}(s^\prime,w)\,\cdot\,[L_v(s^\prime+\frac{1}{2},\overline{F})+L_{v_1}(s^\prime+\frac{1}{2}, \overline{F})\,+$\\

$\dfrac{(q-1)(q^{-w^\prime}-\frac{q^{1-2w^\prime+2s^\prime}}{\alpha\beta})}{(1-\alpha^{-1}q^{1-w^\prime+s^\prime})\,(1-\beta^{-1}q^{1-w^\prime+s^\prime})}\,\cdot\,L_{v_1}(s^\prime+\frac{1}{2}, \overline{F})\,-\,$\\

$q^{-w^\prime}L_{v_1}(s^\prime+w^\prime-\frac{1}{2},\overline{F})]\,\cdot\,F+\frac{1}{4\pi i \kappa} \displaystyle\sum_\chi \int_{\Re(s)=\frac{1}{2}} \Bigl( \displaystyle\int_{Z_\infty \backslash G_\infty} \varphi_\infty \,\cdot\, \overline{W}^E_{s,\chi,\infty} \Bigr)\,\cdot$\\

$\Bigl(\dfrac{L_v(s^\prime+s, \chi)\,\cdot \, L_v(s^\prime+1-s,\overline{\chi})\,\cdot\, |\mathfrak{d}_v|^{-(s^\prime+s-\frac{1}{2})}\,\cdot\, \overline{\chi}(\mathfrak{d}_v)}{L_v(2-2s,\overline{\chi}^2)}\,+\,$\\

$$\dfrac{L_v(s^\prime+s, \chi)\,\cdot \, L_v(s^\prime+1-s,\overline{\chi})\,\cdot\, |\mathfrak{d}_v|^{-(s^\prime+s-\frac{1}{2})}\,\cdot\, \overline{\chi}(\mathfrak{d}_v)}{L_v(2-2s,\overline{\chi}^2)}\,+\,$$\\

\mbox{$\dfrac{L_v(2s-1,\overline{\chi})\,\cdot\,[\dfrac{(q-1)\,q^{1-w^\prime-s+s^\prime}}{(1-q^{2-w^\prime-s+s^\prime})(1-q^{-1+s-s^\prime})}\,-\,\dfrac{q^{-w^\prime}}{1-q^{-w^\prime+s-s^\prime}}]\,\cdot\, |\mathfrak{d}_v|^{\frac{3}{2}-2s}\,\cdot \, \chi(\mathfrak{d}_v)}{L_v(2-2s,\chi^2)}\Bigr)\,\cdot \, E_{s,\chi}\, ds$}\\

\noindent
where
\begin{align*}
\int_{N_v} \varphi_v =
\begin{cases}
\sqrt{\pi}\dfrac{\Gamma(\dfrac{w-1}{2})}{\Gamma(\dfrac{w}{2})} & \text {($v \approx \mathbb{R}$)}\\
2\pi(w-1)^{-1} & \text{($v \approx \mathbb{C}$)}
\end{cases}
\end{align*}
and
\begin{align*}
\int_{Z_v \backslash G_v} \varphi_v \,\cdot\, \overline{W}^E_{s,\chi,v}=
\begin{cases}
\dfrac{\mathcal{G}_v(s,s^\prime,w)}{\pi^{-s}\Gamma(s)}  &\text{($v \approx \mathbb{R}$)}\\\dfrac{\mathcal{G}_v(s,s^\prime,w)}{2 \pi^{-2s-1}\Gamma(2s)}  &\text{($v \approx \mathbb{C}$)}
\end{cases}
\end{align*}

From the spectral decomposition of the Poincar\'e series, and the proof of theorem 4.17 in [Diaconu-Garrett 2008], the Poincar\'e series has meromorphic continuation to a region in $\mathbb{C}^2$ containing $s^\prime=0,\,\,\,w^\prime=1$. As a function of $w^\prime$, for $s^\prime=0$, it is holomorphic in the half-plane $\Re(w^\prime)=\frac{11}{18}$ ([Kim-Shahidi 2002] and [Kim 2005]), except for $w^\prime=1$ where it has a pole of order 1.\\

\vspace{4mm}
\noindent
{\bf 5.1 PRELIMINARIES TO SUBCONVEXITY}\\
\noindent
{\bf 5.1.1 Prologue}\\

Fix a non-archimedean place $v_1$, and take $1<\beta^\prime<2$. Recall that
$$\mathcal{K}_{v_1}(w^\prime, \chi_{v_1}) \ll (q^N)^{-w^\prime}$$
where $\mathcal{K}_{v_1}(w^\prime, \chi_{v_1})$ is the non-decoupled integral for finite prime $v_1$ at which $\chi$ has ramification with conductor $q^N$. Define
$$Z(w^\prime)=\sum_{\chi \in \hat{C}_{0,S}}\, \int_{-\infty}^\infty |L(\frac{1}{2}+it,f \otimes \chi)|^2\,\cdot\,(q^N)^{-w^\prime} \,\cdot \,\mathcal{K}_{\infty}(\frac{1}{2}+it,0,\beta^\prime,\chi)\,dt$$
This is a modified function obtained from
$$I(0,w^\prime)=\sum_\chi\,\frac{1}{2 \pi i}\, \int_{-\infty}^\infty |L(\frac{1}{2}+it,f \otimes \chi)|^2\,\cdot\,\mathcal{K}_{v_1}(w^\prime,\chi_{v_1})\,\cdot \,\mathcal{K}_{\infty}(\frac{1}{2}+it,0,\beta^\prime,\chi)\,dt$$
by taking the asymptotic formula for $\mathcal{K}_{v_1}(w^\prime,\chi_{v_1})$. $Z(w^\prime)$ is absolutely convergent for $\Re(w^\prime)>1$ (see Section 5 in [Diaconu-Garrett 2008]). In this chapter, we will prove the meromorphic continuation and polynomial growth of $Z(w^\prime)$. This will enable us to obtain subconvexity bounds in the $\chi$-depth-aspect.\\

\vspace{4mm}
\noindent
{\bf 5.1.2 Meromorphic continuation of $Z(w^\prime)$}\\

\noindent
{\bf Theorem 5.1}\\
The function
$$Z(w^\prime)=\sum_{\chi \in \hat{C}_{0,S}}\, \int_{-\infty}^\infty |L(\frac{1}{2}+it,f \otimes \chi)|^2\,\cdot\,(q^N)^{-w^\prime} \,\cdot \,\mathcal{K}_{\infty}(\frac{1}{2}+it,0,\beta^\prime,\chi)\,dt$$
where the sum is over a set $\hat{C}_{0,S}$ of characters ramified at the finite prime $v_1$ with conductor $q^N$, and $1<\beta^\prime \leq 2$, $\Re(w^\prime)>1$, has analytic continuation to the half-plane $\Re(w^\prime)>\frac{11}{18}$, except for $w^\prime=1$ where it has a pole of order $1$.\\

\begin{proof}
Let $w^\prime=\delta+i\eta$. Split $Z$ into $Z_1$ and $Z_2$ as follows:
$$Z(w^\prime)=Z_1(w^\prime)+Z_2(w^\prime)$$
Choose a positive constant $C$ and define
$$Z_1(w^\prime)=\sum_{\chi \in \hat{C}_{0,S}:q^N \ll C}\, \int_{-\infty}^\infty |L(\frac{1}{2}+it,f \otimes \chi)|^2\,\cdot\,(q^N)^{-w^\prime} \,\cdot \,\mathcal{K}_{\infty}(\frac{1}{2}+it,0,\beta^\prime,\chi)\,dt$$
We first show that $Z_1(w^\prime)$ has analytic continuation by showing that it is holomorphic for $\delta>0$. Now
$$|Z_1(w^\prime)| \leq \sum_{\chi:q^N \ll C}\, \int_{-\infty}^\infty |L(\frac{1}{2}+it,f \otimes \chi)|^2\,\cdot\,|(q^N)^{-w^\prime}| \,\cdot \,|\mathcal{K}_{\infty}(\frac{1}{2}+it,0,\beta^\prime,\chi)|\,dt$$
$\mathcal{K}_{\infty}(\frac{1}{2}+it,0,\beta^\prime,\chi)$ is positive (see Section 4 in [Diaconu-Garrett 2009]). So
$$|Z_1(w^\prime)| \leq \sum_{\chi:q^N \ll C}\, \int_{-\infty}^\infty |L(\frac{1}{2}+it,f \otimes \chi)|^2\,\cdot\,(q^N)^{-\delta} \,\cdot \,\mathcal{K}_{\infty}(\frac{1}{2}+it,0,\beta^\prime,\chi)\,dt$$
Since
$$(q^N)^{-\delta} \ll_{\beta^\prime,C} (q^N)^{-\beta^\prime}$$
then\\

$|Z_1(w^\prime)| \ll \displaystyle\sum_{\chi:q^N \ll C}\, \displaystyle\int_{-\infty}^\infty |L(\frac{1}{2}+it,f \otimes \chi)|^2\,\cdot\,(q^N)^{-\beta^\prime} \,\cdot \,\mathcal{K}_{\infty}(\frac{1}{2}+it,0,\beta^\prime,\chi)\,dt$\\

$< \displaystyle\sum_{\chi \in \hat{C}_{0,S}}\, \displaystyle\int_{-\infty}^\infty |L(\frac{1}{2}+it,f \otimes \chi)|^2\,\cdot\,(q^N)^{-\beta^\prime} \,\cdot \,\mathcal{K}_{\infty}(\frac{1}{2}+it,0,\beta^\prime,\chi)\,dt$\\

$=Z(\beta^\prime)$\\

\noindent
which is convergent for $\Re(w^\prime)<\frac{2}{9}$. Thus, $Z_1(w^\prime)$ is holomorphic for \mbox{$\Re(w^\prime)=\delta>0$} (in particular for $\Re(w^\prime)>\frac{11}{18}$).\\

Now we prove that $Z_2(w^\prime)$ has analytic continuation. Consider
$$I(s^\prime,w^\prime,\beta^\prime)=\sum_{\chi \in \hat{C}_{0,S}}\, \frac{1}{2\pi i}\,\int_{-\infty}^\infty L(s,f \otimes \chi)\,\cdot\,L(s^\prime+1-s,\overline{f} \otimes \overline{\chi})\, \cdot \, \mathcal{K}_{v_1}(w^\prime,\chi_{v_1})\,\cdot \, \mathcal{K}_{\infty}(s,s^\prime,\beta^\prime,\chi)\,dt$$
where $$\mathcal{K}_{\infty}(s,s^\prime,\beta^\prime,\chi)=\prod_{v|\infty} \mathcal{K}_v(s,s^\prime,\beta^\prime,\chi_v)$$
Recall that this expression is obtained from the integral representation 
$$\int_{Z_\mathbb{A} G_k \backslash G_\mathbb{A}} P \acute{e}\, |f|^2\,dg$$
where $P \acute{e}(g)$ converges absolutely and locally uniformly for $\Re(s^\prime)>1$ and $\Re(w^\prime)>1$. Also recall from the spectral decomposition of $P \acute{e}$, that $P \acute{e}$ has meromorphic continuation to $\Re(w^\prime)>\frac{11}{18}$ with a pole of order $1$ at $w^\prime=1$. Thus
$$I(0,w^\prime,\beta^\prime)=\sum_{\chi \in \hat{C}_{0,S}}\, \frac{1}{2\pi i}\,\int_{-\infty}^\infty |L(\frac{1}{2}+it,f \otimes \chi)|^2\,\cdot\, \mathcal{K}_{v_1}(w^\prime,\chi_{v_1})\,\cdot \, \mathcal{K}_{\infty}(\frac{1}{2}+it,0,\beta^\prime,\chi)\,dt$$
is holomorphic for $\Re(w^\prime)>\frac{11}{18}$ except at $w^\prime=1$ where there is a pole of order $1$.\\
\indent
In the region of absolute convergence for $\Re(w^\prime)=\delta>1$, write
$$I(0,w^\prime,\beta^\prime)\,=\,I_1(0,w^\prime,\beta^\prime)+I_2(0,w^\prime,\beta^\prime)$$
where
$$I_1(0,w^\prime,\beta^\prime)=\sum_{\chi:q^N \ll C}\, \textstyle\frac{1}{2\pi i}\,\displaystyle\int_{-\infty}^\infty |L(\frac{1}{2}+it,f \otimes \chi)|^2\,\cdot\, \mathcal{K}_{v_1}(w^\prime,\chi_{v_1})\,\cdot \, \mathcal{K}_{\infty}(\textstyle\frac{1}{2}+it,0,\beta^\prime,\chi)\,dt$$
Now
$$I(0,w^\prime,\beta^\prime)\,=\,I_1(0,w^\prime,\beta^\prime)\,+\sum_{\chi:q^N \gg C}\, \textstyle\frac{1}{2\pi i}\,\int_{-\infty}^\infty |L(\frac{1}{2}+it,f \otimes \chi)|^2\,\cdot\, C^\prime \, \cdot \, (q^N)^{-w^\prime}\,\cdot \, \mathcal{K}_{\infty}(\textstyle\frac{1}{2}+it,0,\beta^\prime,\chi)\,dt$$
where the constant
$$C^\prime=\frac{q}{q-1}\, \cdot \, \frac{1-|\alpha|^2 |\beta|^2 q^{-2w^\prime}}{(1-|\alpha|^2q^{-w^\prime})(1-|\beta|^2q^{-w^\prime})(1-\overline{\alpha}\beta q^{-w^\prime})(1-\alpha \overline{\beta}q^{-w^\prime})}$$
since
$$\mathcal{K}_{v_1}(w^\prime,\chi_{v_1})=\frac{q^{1-Nw^\prime}}{q-1}\, \cdot \, \frac{1-|\alpha|^2 |\beta|^2 q^{-2w^\prime}}{(1-|\alpha|^2q^{-w^\prime})(1-|\beta|^2q^{-w^\prime})(1-\overline{\alpha}\beta q^{-w^\prime})(1-\alpha \overline{\beta}q^{-w^\prime})}$$
So
$$I(0,w^\prime,\beta^\prime)=I_1(0,w^\prime,\beta^\prime)+C^\prime \cdot Z_2(w^\prime)$$
Thus, to show that $Z_2(w^\prime)$ has analytic continuation, it suffices to show that $I_1(0,w^\prime,\beta^\prime)$ is absolutely convergent for $\Re(w^\prime)>\frac{11}{18}$.\\

$I_1(0,w^\prime,\beta^\prime)=$\\

$\displaystyle\sum_{\chi:q^N \ll C}\, \frac{1}{2\pi i}\,\displaystyle\int_{-\infty}^\infty |L(\frac{1}{2}+it,f \otimes \chi)|^2\,\cdot\, \mathcal{K}_{v_1}(w^\prime,\chi_{v_1})\,\cdot \, \mathcal{K}_{\infty}(\frac{1}{2}+it,0,\beta^\prime,\chi)\,dt$\\

$\ll \displaystyle\sum_{\chi:q^N \ll C}\, \frac{1}{2\pi i}\,\displaystyle\int_{-\infty}^\infty |L(\frac{1}{2}+it,f \otimes \chi)|^2\,\cdot\, |(q^N)^{-w^\prime}|\,\cdot \, \mathcal{K}_{\infty}(\frac{1}{2}+it,0,\beta^\prime,\chi)\,dt$\\

$\ll \displaystyle\sum_{\chi:q^N \ll C}\, \frac{1}{2\pi i}\,\displaystyle\int_{-\infty}^\infty |L(\frac{1}{2}+it,f \otimes \chi)|^2\,\cdot\, (q^N)^{-\beta}\,\cdot \, \mathcal{K}_{\infty}(\frac{1}{2}+it,0,\beta^\prime,\chi)\,dt$\\

$< \displaystyle\sum_{\chi \in \hat{C}_{0,S}}\, \frac{1}{2\pi i}\,\displaystyle\int_{-\infty}^\infty |L(\frac{1}{2}+it,f \otimes \chi)|^2\,\cdot\, (q^N)^{-\beta}\,\cdot \, \mathcal{K}_{\infty}(\frac{1}{2}+it,0,\beta^\prime,\chi)\,dt$\\

$=Z(\beta^\prime)$\\

\noindent
which is convergent for $\Re(w^\prime)>\frac{2}{9}$. Thus $Z_2(w^\prime)$ is absolutely convergent for $\Re(w^\prime)>\frac{11}{18}$, proving the theorem.\\
\end{proof}

\vspace{4mm}
\newpage
\noindent
{\bf 5.1.3 Polynomial growth of $Z(w^\prime)$}\\

\noindent
{\bf Theorem 5.2}\\
For every fixed small positive $\epsilon$, the generating function
$$Z(w^\prime)=\sum_{\chi \in \hat{C}_{0,S}}\, \int_{-\infty}^\infty |L(\frac{1}{2}+it,f \otimes \chi)|^2\,\cdot\,(q^N)^{-w^\prime} \,\cdot \,\mathcal{K}_{\infty}(\frac{1}{2}+it,0,\beta^\prime,\chi)\,dt$$
has polynomial growth in the conductor $q^N$ for $\frac{11}{18}+\epsilon \leq \Re(w^\prime) \leq 1+\epsilon$; that is, on the vertical line $\Re(w^\prime)=\frac{11}{18}+\epsilon$,
$$Z(w^\prime) \ll_{\epsilon, \beta^\prime} (q^N)^\gamma$$
with a computable $\gamma>0$ independent of $\beta^\prime$.\\

\noindent
Define the Poincar\'e series data at the non-archimedean place $v=v_1$ as earlier, namely:
\begin{align*}
\varphi(\begin{smallmatrix} 1 & x\\0 & 1 \end{smallmatrix})=
\begin{cases}
1\,\,\,\,\, & \text {($x \in \mathfrak{o}_v$)}\\
|x|^{-w^\prime}\,\,\,\,\, & \text{($x \not \in \mathfrak{o}_v$)}
\end{cases}
\end{align*}

\noindent
Denote
$$I(s^\prime,w^\prime,\beta^\prime)=\int_{Z_\mathbb{A} G_k \backslash G_\mathbb{A}} P \acute{e}(g)\, |f(g)|^2\,dg$$
For $\Re(s^\prime),\Re(w^\prime)>1$,
$$I(s^\prime,w^\prime,\beta^\prime)=\sum_{\chi \in \hat{C}_{0,S}}\, \frac{1}{2\pi i}\,\int_{-\infty}^\infty L(s,f \otimes \chi)\,\cdot\,L(s^\prime+1-s,\overline{f} \otimes \overline{\chi})\, \cdot \, \mathcal{K}_{v_1}(w^\prime,\chi_{v_1})\,\cdot \, \mathcal{K}_{\infty}(s,s^\prime,\beta^\prime,\chi)\,dt$$
where $$\mathcal{K}_{\infty}(s,s^\prime,\beta^\prime,\chi)=\prod_{v|\infty} \mathcal{K}_v(s,s^\prime,\beta^\prime,\chi_v)$$
In the region of absolute convergence
$$I(0,w^\prime,\beta^\prime)=I_1(0,w^\prime,\beta^\prime)+I_2(0,w^\prime,\beta^\prime)$$
where
$$I_1(0,w^\prime,\beta^\prime)\,=\,\sum_{\chi:q^N \ll C}\, \frac{1}{2\pi i}\,\int_{-\infty}^\infty |L(\frac{1}{2}+it,f \otimes \chi)|^2\,\cdot\, \mathcal{K}_{v_1}(w^\prime,\chi_{v_1})\,\cdot \, \mathcal{K}_{\infty}(\frac{1}{2}+it,0,\beta^\prime,\chi)\,dt$$
and 
$$I_2(0,w^\prime,\beta^\prime)=C^\prime \, \cdot \, Z_2(w^\prime)$$
The constant
$$C^\prime=\frac{q}{q-1}\, \cdot \, \frac{1-|\alpha|^2 |\beta|^2 q^{-2w^\prime}}{(1-|\alpha|^2 q^{-w^\prime})(1-|\beta|^2 q^{-w^\prime})(1-\overline{\alpha}\beta q^{-w^\prime})(1-\alpha \overline{\beta}q^{-w^\prime})}$$
and
$$Z_2(w^\prime)=\sum_{\chi:q^N \gg C}\, \int_{-\infty}^\infty|L(\frac{1}{2}+it,f \otimes \chi)|^2\,\cdot\,(q^N)^{-w^\prime} \,\cdot \,\mathcal{K}_{\infty}(\frac{1}{2}+it,0,\beta^\prime,\chi)\,dt$$
with
$$Z_2(w^\prime)=Z(w^\prime)-Z_1(w^\prime)$$
so
$$I_2(0,w^\prime,\beta^\prime)=C^\prime[Z(w^\prime)-Z_1(w^\prime)]$$
Recall that $Z_1(w^\prime)$ is holomorphic in the half-plane $\Re(w^\prime)>\frac{11}{18}$, and that $C^\prime$ is the positive constant where the cutoff of $I(0,w^\prime,\beta^\prime)$ was made. So $Z_1(w^\prime)$ has polynomial growth in $q^N$. Thus, the polynomial bound of $Z(w^\prime)$ will be deduced from that of $I_2(0,w^\prime,\beta^\prime)$. That is, we will prove that $I_2(0,w^\prime,\beta^\prime)$ has polynomial growth in $q^N$.\\

\noindent
Again recall that
$$I(w^\prime)=I(0,w^\prime,\beta^\prime)=\sum_{\chi \in \hat{C}_{0,S}}\, \frac{1}{2\pi i}\,\int_{-\infty}^\infty |L(\frac{1}{2}+it,f \otimes \chi)|^2\,\cdot\, \mathcal{K}_{v_1}(w^\prime,\chi_{v_1})\,\cdot \, \mathcal{K}_{\infty}(\frac{1}{2}+it,0,\beta^\prime,\chi)\,dt$$ 
which is obtained from the integral representation
$$\int_{Z_\mathbb{A} G_k \backslash G_\mathbb{A}} P \acute{e}\, |f|^2\,dg$$
$P \acute{e}$ admits a spectral decomposition
$$P \acute{e}=\,\,\,\mbox{Singular part}\,\,\,+\,\,\,\mbox{Cuspidal part}\,\,\,+\,\,\,\mbox{Continuous part}$$
In the spectral decomposition set $s^\prime=0$ and obtain\\

$P\acute{e}=$\\
$\lim_{s^\prime \to 0}\Bigl(\displaystyle\int_{N_\infty} \varphi_\infty \Bigr)\,\cdot\,\dfrac{1-q^{-w^\prime}}{1-q^{1-w^\prime}}  \,\cdot\, E_{s^\prime+1,1}+\displaystyle\sum_F \overline{\rho}_F \mathcal{G}_{F_\infty}(\beta^\prime)\,\cdot\,\Bigl[L_v(\frac{1}{2},\overline{F})+$\\

\mbox{$L_{v_1}(\frac{1}{2}, \overline{F})\,+\,\dfrac{(q-1)(q^{-w^\prime}-\frac{q^{1-2w^\prime}}{\alpha\beta})}{(1-\alpha^{-1}q^{1-w^\prime})\,(1-\beta^{-1}q^{1-w^\prime})}\,\cdot\,L_{v_1}(\frac{1}{2}, \overline{F})\,-\,q^{-w^\prime}L_{v_1}(\frac{2w^\prime-1}{2},\overline{F})\Bigr]\,\cdot\,F+$}\\

$\displaystyle\sum_\chi \dfrac{\overline{\chi}(\mathfrak{d}_v)}{4 \pi i \kappa} \displaystyle\int_{\Re(s)=\frac{1}{2}} \Bigl( \int_{Z_\infty \backslash G_\infty} \varphi_\infty \,\cdot\, \overline{W}^E_{s,\overline{\chi},\infty} \Bigr)\,\cdot$\\

$\Bigl(\dfrac{L_v(s, \chi)\,\cdot \, L_v(1-s,\overline{\chi})\,\cdot\, |\mathfrak{d}_v|^{-(s-\frac{1}{2})}}{L_v(2-2s,\overline{\chi}^2)}\,+\,\dfrac{L_{v_1}(s,\chi)\,\cdot \, L_{v_1}(1-s,\overline{\chi})\,\cdot\, |\mathfrak{d}_{v_1}|^{-(s-\frac{1}{2})}}{L_{v_1}(2-2s,\overline{\chi}^2)}\,+\,$\\

\mbox{$\dfrac{L_{v_1}(2s-1,\overline{\chi})\,\cdot\,\Bigl[\dfrac{(q-1)\,q^{1-w^\prime-s}}{(1-q^{2-w^\prime-s})(1-q^{s-1})}\,-\,\dfrac{q^{-w^\prime}}{1-q^{-w^\prime+s}}\Bigr]\,\cdot\, |\mathfrak{d}_{v_1}|^{\frac{3}{2}-2s}}{L_{v_1}(2-2s,\chi^2)}\Bigr)\,\cdot \, E_{s,\chi}\, ds$}\\

\noindent
So
\newcommand{\sing}{sing}
\newcommand{\cont}{cont}
$$I(w^\prime)=I_{sing}(w^\prime)+I_{cusp}(w^\prime)+I_{cont}(w^\prime)$$
where\\

$I_{sing}(w^\prime)=\lim_{s^\prime \to 0}\Bigl(\displaystyle\int_{N_\infty} \varphi_\infty \Bigr) \, \cdot \, \dfrac{1-q^{-w^\prime}}{1-q^{1-w^\prime}} \,\cdot\, \langle E_{s^\prime+1,1},|f|^2 \rangle$\\

$I_{cusp}(w^\prime)=\displaystyle\sum_F \overline{\rho}_F \mathcal{G}_{F_\infty}(\beta^\prime)\,\cdot\,\Bigl[2L_v(\frac{1}{2},\overline{F})+\dfrac{(q-1)(q^{-w^\prime}-\frac{q^{1-2w^\prime}}{\alpha\beta})}{(1-\alpha^{-1}q^{1-w^\prime})\,(1-\beta^{-1}q^{1-w^\prime})}\,\cdot\,$\\

$$L_{v}(\frac{1}{2}, \overline{F})\,-\,q^{-w^\prime}L_{v}(\frac{2w^\prime-1}{2},\overline{F})\Bigr]\,\cdot\,\langle F,|f|^2 \rangle$$

$I_{cont}(w^\prime)=\displaystyle\sum_\chi \dfrac{\overline{\chi}(\mathfrak{d}_v)}{4\pi i \kappa} \displaystyle\int_{\Re(s)=\frac{1}{2}} \Bigl( \displaystyle\int_{Z_\infty \backslash G_\infty} \varphi_\infty \,\cdot\, \overline{W}^E_{1-s,\overline{\chi},\infty} \Bigr)\,\cdot\,\Bigl(\dfrac{2L_v(s, \chi)\,\cdot \, L_v(1-s,\overline{\chi})\,\cdot\, |\mathfrak{d}_v|^{-(s-\frac{1}{2})}}{L_v(2-2s,\overline{\chi}^2)}\,+\,$\\

$\dfrac{L_{v}(2s-1,\overline{\chi})\,\cdot\,\Bigl[\dfrac{(q-1)\,q^{1-w^\prime-s}}{(1-q^{2-w^\prime-s})(1-q^{s-1})}\,-\,\dfrac{q^{-w^\prime}}{1-q^{-w^\prime+s}}\Bigr]\,\cdot\, |\mathfrak{d}_{v}|^{\frac{3}{2}-2s}}{L_{v}(2-2s,\chi^2)}\Bigr)\,\cdot \,\langle E_{s,\chi},|f|^2 \rangle \, ds$\\

\noindent
Note that the dependence on $w^\prime$ is at $v_1$ only. Let

\newcommand{\aux}{aux}
\noindent
$$\mathcal{M}_1(w^\prime)=\dfrac{(q-1)(q^{-w^\prime}-\frac{q^{1-2w^\prime}}{\alpha\beta})}{(1-\alpha^{-1}q^{1-w^\prime})\,(1-\beta^{-1}q^{1-w^\prime})}\,\cdot\,L(\frac{1}{2}, \overline{F})\,-\,q^{-w^\prime}L(\frac{2w^\prime-1}{2},\overline{F})$$
\noindent
and\\

$\mathcal{M}_2(w^\prime)=\dfrac{(q-1)\,q^{1-w^\prime-s}}{(1-q^{2-w^\prime-s})(1-q^{s-1})}\,-\,\dfrac{q^{-w^\prime}}{1-q^{-w^\prime+s}}$\\

\noindent
Then define the auxiliary function $I^{\aux}(w^\prime)$ by\\

$I^{\aux}(w^\prime)=\displaystyle\sum_F \overline{\rho}_F \mathcal{G}_{F_\infty}(\beta^\prime)\,\cdot\,[2L(\frac{1}{2},\overline{F})+\mathcal{M}_1^{\aux}(w^\prime)]\,\cdot\,\langle F,|f|^2 \rangle\, + \,$\\

$\displaystyle\sum_\chi \dfrac{\overline{\chi}(\mathfrak{d})}{4\pi i \kappa} \displaystyle\int_{\Re(s)=\frac{1}{2}} \Bigl( \int_{Z_\infty \backslash G_\infty} \varphi_\infty \,\cdot\, \overline{W}^E_{1-s,\overline{\chi},\infty} \Bigr)\,\cdot\,\dfrac{2L(s,\chi)\,\cdot \, L(1-s,\overline{\chi})\,\cdot\, |\mathfrak{d}|^{-(s-\frac{1}{2})}}{L(2-2s,\overline{\chi}^2)}\,+\,$\\

$\dfrac{L(2s-1,\overline{\chi})\,\cdot\,\mathcal{M}_2^{\aux}(w^\prime)\,\cdot\, |\mathfrak{d}|^{\frac{3}{2}-2s}}{L(2-2s,\chi^2)}\,\cdot \,\langle E_{s,\chi},|f|^2 \rangle \, ds$\\

\noindent
where $\mathcal{M}_1^{\aux}(w^\prime)$ and $\mathcal{M}_2^{\aux}(w^\prime)$ are defined by
$$\mathcal{M}_1^{\aux}(w^\prime)=\mathcal{M}_1(w^\prime)\, \cdot \, (q^N)^\gamma$$
$$\mathcal{M}_2^{\aux}(w^\prime)=\mathcal{M}_2(w^\prime)\, \cdot \, (q^N)^\gamma$$
where $\gamma>0$, independent of $\beta^\prime$. Define\\

$H(w^\prime)=I(w^\prime)-I^{\aux}(w^\prime)\,=\,\lim_{s^\prime \to 0}\Bigl(\displaystyle\int_{N_\infty} \varphi_\infty \Bigr) \, \cdot \, \dfrac{1-q^{-w^\prime}}{1-q^{1-w^\prime}} \,\cdot\, \langle E_{s^\prime+1,1},|f|^2 \rangle\, + \,$\\

$\displaystyle\sum_F \overline{\rho}_F \mathcal{G}_{F_\infty}(\beta^\prime)\,\cdot\,[\mathcal{M}_1(w^\prime)\,-\, \mathcal{M}_1^{\aux}(w^\prime)]\,\cdot\,\langle F,|f|^2 \rangle\, + \,\mbox{$\displaystyle\sum_\chi \dfrac{\overline{\chi}(\mathfrak{d})}{4\pi i \kappa} \int_{\Re(s)=\frac{1}{2}} \Bigl( \int_{Z_\infty \backslash G_\infty} \varphi_\infty \,\cdot\, \overline{W}^E_{1-s,\overline{\chi},\infty} \Bigr)\,\cdot\,$}$\\

$\dfrac{L(2s-1,\overline{\chi})\,\cdot\,[\mathcal{M}_2(w^\prime)\, - \,\mathcal{M}_2^{\aux}(w^\prime)]\,\cdot\, |\mathfrak{d}|^{\frac{3}{2}-2s}}{L(2-2s,\chi^2)}\,\cdot \,\langle E_{s,\chi},|f|^2 \rangle \, ds$\\
\\

\noindent
{\bf Proposition 5.3}\\
For $\epsilon$ sufficiently small, $$H(w^\prime)=I(w^\prime)-I^{\aux}(w^\prime)$$ restricted to $\frac{11}{18} < \Re(w^\prime) \leq 1+\epsilon$, extends holomorphically to the whole vertical strip \mbox{$-\epsilon \leq \Re(w^\prime) \leq 1+\epsilon$.}\\

\begin{proof}

The first term in $H(w^\prime)$, i.e. $$\lim_{s^\prime \to 0}\Bigl(\int_{N_\infty} \varphi_\infty \Bigr) \, \cdot \, \dfrac{1-q^{-w^\prime}}{1-q^{1-w^\prime}} \,\cdot\, \langle E_{s^\prime+1,1},|f|^2 \rangle$$ is holomorphic in the strip $-\epsilon \leq \Re(w^\prime) \leq 1+\epsilon$, except at $w^\prime=0,1$ where there are poles.\\  
$\mathcal{M}_1(w^\prime)-\mathcal{M}_1^{\aux}(w^\prime)=\mathcal{M}_1(w^\prime)-\mathcal{M}_1(w^\prime)\, \cdot \, (q^N)^\gamma=\mathcal{M}_1(w^\prime)[1-(q^N)^\gamma]$ and \\
\noindent
$\mathcal{M}_2(w^\prime)-\mathcal{M}_2^{\aux}(w^\prime)=\mathcal{M}_2(w^\prime)-\mathcal{M}_2(w^\prime)\, \cdot \, (q^N)^\gamma=\mathcal{M}_2(w^\prime)[1-(q^N)^\gamma]$\\
\noindent
Since both $\mathcal{M}_1(w^\prime)$ and $\mathcal{M}_2(w^\prime)$ are holomorphic in the strip, then $H(w^\prime)$ is also holomorphic in the strip.\\
\end{proof}

\noindent
{\bf Proposition 5.4}\\
Fix a small positive $\epsilon$. For $\frac{11}{18}+\epsilon \leq \Re(w^\prime) \leq 1+\epsilon$, or $\Re(w^\prime)=-\epsilon$,
$$I^{\aux}(w^\prime) \ll_{\epsilon,\beta^\prime} (q^N)^\gamma$$

\begin{proof}

Again,\\

$I^{\aux}(w^\prime)\,=\,\displaystyle\sum_F \overline{\rho}_F \mathcal{G}_{F_\infty}(\beta^\prime)\,\cdot\,[2L(\frac{1}{2},\overline{F})+\mathcal{M}_1^{\aux}(w^\prime)\,\cdot\,\langle F,|f|^2 \rangle\, + \,$\\

$\displaystyle\sum_\chi \dfrac{\overline{\chi}(\mathfrak{d})}{4\pi i \kappa} \displaystyle\int_{\Re(s)=\frac{1}{2}} \Bigl( \displaystyle\int_{Z_\infty \backslash G_\infty} \varphi_\infty \,\cdot\, \overline{W}^E_{1-s,\overline{\chi},\infty} \Bigr)\,\cdot\,\dfrac{2L(s,\chi)\,\cdot \, L(1-s,\overline{\chi})\,\cdot\, |\mathfrak{d}|^{-(s-\frac{1}{2})}}{L(2-2s,\overline{\chi}^2)}\,+\,$\\

$\dfrac{L(2s-1,\overline{\chi})\,\cdot\,\mathcal{M}_2^{\aux}(w^\prime)\,\cdot\, |\mathfrak{d}|^{\frac{3}{2}-2s}}{L(2-2s,\chi^2)}\,\cdot \,\langle E_{s,\chi},|f|^2 \rangle \, ds$\\

\noindent
where\\
$$\mathcal{M}_1^{\aux}(w^\prime)=\mathcal{M}_1(w^\prime)\, \cdot \, (q^N)^\gamma$$
$$\mathcal{M}_2^{\aux}(w^\prime)=\mathcal{M}_2(w^\prime)\, \cdot \, (q^N)^\gamma$$

$\mathcal{M}_1(w^\prime)=\dfrac{(q-1)(q^{-w^\prime}-\frac{q^{1-2w^\prime}}{\alpha\beta})}{(1-\alpha^{-1}q^{1-w^\prime})\,(1-\beta^{-1}q^{1-w^\prime})}\,\cdot\,L(\frac{1}{2}, \overline{F})\,-\,q^{-w^\prime}L(\frac{2w^\prime-1}{2},\overline{F})$\\

$\mathcal{M}_2(w^\prime)=\dfrac{(q-1)\,q^{1-w^\prime-s}}{(1-q^{2-w^\prime-s})(1-q^{s-1})}\,-\,\dfrac{q^{-w^\prime}}{1-q^{-w^\prime+s}}$\\

$\mathcal{M}_1^{\aux}(w^\prime) \ll (q^N)^\gamma$\\

$\mathcal{M}_2^{\aux}(w^\prime) \ll (q^N)^\gamma$\\

All other terms in $I^{\aux}(w^\prime)$ are independent of the conductor $q^N$, and have a polynomial bound. Thus
$$I^{\aux}(w^\prime) \ll_{\epsilon,\beta^\prime} (q^N)^\gamma$$

\end{proof}

Recall we are trying to prove a polynomial bound for $I_2(w^\prime)$ in the conductor $q^N$. Now 
$$I_2(w^\prime)=I_2(w^\prime)-I^{\aux}(w^\prime)+I^{\aux}(w^\prime)$$
We have proven a polynomial bound for $I^{\aux}(w^\prime)$, so we now prove a polynomial bound for $I_2(w^\prime)-I^{\aux}(w^\prime)$.\\
$$H(w^\prime)-I_1(w^\prime)=I_2(w^\prime)-I^{\aux}(w^\prime)$$
Thus it suffices to prove a polynomial bound for $H(w^\prime)-I_1(w^\prime)$ on the line $\Re(w^\prime)=\frac{11}{18}+\epsilon$.\\

\begin{proof}

Recall that $H(w^\prime)$ is holomorphic in the strip $-\epsilon < \Re(w^\prime) <1+\epsilon$. Also recall
$$I_1(w^\prime)=\sum_{\chi:q^N \ll C}\, \frac{1}{2\pi i}\,\int_{-\infty}^\infty |L(\frac{1}{2}+it,f \otimes \chi)|^2\,\cdot\, \mathcal{K}_{v_1}(w^\prime,\chi_{v_1})\,\cdot \, \mathcal{K}_{\infty}(\frac{1}{2}+it,0,\beta^\prime,\chi)\,dt$$
and
$$I_1(w^\prime) \ll Z(\beta^\prime) < \infty$$
So $I_1(w^\prime)$ converges absolutely throughout the strip; i.e. $I_1(w^\prime)$ is holomorphic throughout the strip. Thus $H(w^\prime)-I_1(w^\prime)$ is also holomorphic throughout the strip. Recall
$$I_2(w^\prime)=C^\prime[Z(w^\prime)-Z_1(w^\prime)]$$
For $\Re(w^\prime)=1+\epsilon$, since $I^{\aux}(w^\prime) \ll (q^N)^\gamma$, for $\gamma>0$, $Z(w^\prime)=O(1)$ and $Z_1(w^\prime)$ already has polynomial growth in $q^N$, we conclude that $$H(w^\prime)-I_1(w^\prime)=I_2(w^\prime)-I^{\aux}(w^\prime)$$ has polynomial growth in $q^N$ for $\Re(w^\prime)=1+\epsilon$.\\
\indent
Now assume $\Re(w^\prime)=-\epsilon$.
$$H(w^\prime)-I_1(w^\prime)=I(w^\prime)-I^{\aux}(w^\prime)-I_1(w^\prime)$$
Again, $I^{\aux}(w^\prime)$ has polynomial growth for $\Re(w^\prime)=-\epsilon$, and $I_1(w^\prime) \ll Z(\beta^\prime)$. The spectral expansion of $I(w^\prime)$ and $I_1(w^\prime)$ shows that $I(w^\prime)$ and $I_1(w^\prime)$ also have polynomial growth for $\Re(w^\prime)=-\epsilon$. Thus $H(w^\prime)-I_1(w^\prime)$ has polynomial growth in $q^N$ for \mbox{$\Re(w^\prime)=-\epsilon$.}\\
\end{proof}

We now apply Phragmen-Lindel\"of and conclude that $$I_2(w^\prime)-I^{\aux}(w^\prime)$$ has polynomial growth in $q^N$ within the strip $\frac{11}{18}+\epsilon \leq \Re(w^\prime) \leq 1+\epsilon$, and hence, so has $I_2(w^\prime)$. Thus, we have proven that $Z(w^\prime)$ has polynomial growth in $q^N$ within the strip $\frac{11}{18}+\epsilon \leq \Re(w^\prime) \leq 1+\epsilon.$\\

\vspace{4mm}
\noindent
{\bf 6.1 SUBCONVEXITY BOUNDS}\\
\noindent
{\bf 6.1.1 Prologue}\\

Our goal is to break convexity in the $\chi$-depth-aspect for a family of $L$-functions $L(\frac{1}{2}+it, f \otimes \chi)$, where $\chi$ has arbitrary ramification at a fixed finite prime $v_1$. For a cuspform $f$ on $GL_2(k)$, the $\chi$-depth-aspect convexity bound for the twisted $L$-function $L(\frac{1}{2}+it,f \otimes \chi)$ is
$$L(\frac{1}{2}+it,f \otimes \chi) \ll q^{N(\frac{d}{2}+\epsilon)}$$
where $q^N$ with $N \geq 1$ is the conductor of $\chi$ allowed to be ramified at the finite place $v_1$, and $d$ is the degree of the number field $k$ over $\mathbb{Q}$. Using methods in section 4 in [Diaconu-Garrett 2009], we will break convexity at the finite place $v_1$ by decreasing the exponent. So fix a non-archimedean place $v_1$, take $1<\beta^\prime<2$ and fix $0<t<1$ in the nondecoupled integral at the archimedean places.\\
\noindent
Write
$$I(0,w^\prime)=\sum_{\chi \in \hat{C}_{0,S}}\, \frac{1}{2\pi i}\,\int_{-\infty}^\infty |L(\frac{1}{2}+it,f \otimes \chi)|^2\,\cdot\, \mathcal{K}_{v_1}(w^\prime,\chi_{v_1})\,\cdot \, \mathcal{K}_{\infty}(\frac{1}{2}+it,0,\beta^\prime,\chi)\,dt$$
where $$\mathcal{K}_{\infty}(\frac{1}{2}+it,0,\beta^\prime,\chi)=\prod_{v|\infty} \mathcal{K}_v(\frac{1}{2}+it,0,\beta^\prime,\chi_v)$$
is the nondecoupled integral at the archimedean places, and $\mathcal{K}_{v_1}(w^\prime,\chi_{v_1})$ is the nondecoupled integral at the finite prime $v_1$; $\mathcal{K}_{v_1}(w^\prime,\chi_{v_1})$ does not depend on $t$. We have shown that
$$\mathcal{K}_{v_1}(w^\prime,\chi_{v_1})=\frac{q^{1-Nw^\prime}}{q-1}\, \cdot \, \frac{1-|\alpha|^2 |\beta|^2 q^{-2w^\prime}}{(1-|\alpha|^2q^{-w^\prime})(1-|\beta|^2q^{-w^\prime})(1-\overline{\alpha}\beta q^{-w^\prime})(1-\alpha \overline{\beta}q^{-w^\prime})}$$
and
$$\mathcal{K}_{v_1}(w^\prime,\chi_{v_1}) \ll (q^N)^{-w^\prime}$$
Define
$$Z(w^\prime)=\sum_{\chi \in \hat{C}_{0,S}}\, \int_{-\infty}^\infty |L(\frac{1}{2}+it,f \otimes \chi)|^2\,\cdot\,(q^N)^{-w^\prime} \,\cdot \,\mathcal{K}_{\infty}(\frac{1}{2}+it,0,\beta^\prime,\chi)\,dt$$
$Z(w^\prime)$ has analytic continuation to $\Re(w^\prime)>\frac{11}{18}$ with a pole of order $1$ at \mbox{$w^\prime=1$}, and has polynomial growth on every vertical strip inside \mbox{$\frac{11}{18}+\epsilon \leq \Re(w^\prime) \leq 1+\epsilon$}. Choose $\frac{11}{18}<\delta_0<1$. From section 4 in [Diaconu-Garrett 2008], for $\delta_0 \leq \Re(w^\prime) \leq 1+\epsilon$, by Phragmen-Lindel\"of, $Z(\delta_0+i\eta)$ has polynomial growth of exponent less than $\frac{1}{2}$. Consider the rectangle R with vertices at $\delta_0-iS\,$, \mbox{$\beta^\prime-iS,$} $\,\beta^\prime+iS, \, \delta_0+iS$. Recall Perron's formula: for $\beta^\prime>1$,
\begin{align*}
\frac{1}{2\pi i}\int_{\beta^\prime-iS}^{\beta^\prime+iS}\frac{x^w}{w}\,dw=
\begin{cases}
1 & \text{(for $x>1$)}\\
0 & \text{(for $x<1$)}
\end{cases}
\,+\, x^{\beta^\prime} O_{\beta^\prime}(\min\{1,\frac{1}{S|\log x|}\})
\end{align*}
Applying Perron's formula to the integral
$$\frac{1}{2\pi i}\int_{\beta^\prime-iS}^{\beta^\prime+iS}\frac{Z(w^\prime)\,x^{w^\prime}}{w^\prime}\,dw^\prime$$
gives\\

$\frac{1}{2\pi i}\displaystyle \int_{\beta^\prime-iS}^{\beta^\prime+iS}\dfrac{Z(w^\prime)\,x^{w^\prime}}{w^\prime}\,dw^\prime$\\

$=\frac{1}{2\pi i}\,\displaystyle \sum_\chi\,\displaystyle \int_{-\infty}^\infty |L(\frac{1}{2}+it,f \otimes \chi)|^2\,\cdot\, \Bigl(\displaystyle \int_{\beta-iS}^{\beta+iS}\dfrac{(x/q^N)^{w^\prime}}{w^\prime}\,dw^\prime \Bigr)\,\cdot \,$\\
$$\mathcal{K}_{\infty}(\frac{1}{2}+it,0,\beta^\prime,\chi)\,dt$$

$=\displaystyle \sum_{\chi: q^N\leq x}\,\displaystyle \int_{-\infty}^\infty \, |L(\frac{1}{2}+it,f \otimes \chi)|^2\,\cdot\,1\,\cdot\,\mathcal{K}_{\infty}(\frac{1}{2}+it,0,\beta^\prime,\chi)\,dt\, + \, E(x,S)$\\

\noindent
where the error term $E(x,S)$ is 
$$E(x,S) \ll \sum_{\chi}\,\int_{-\infty}^\infty |L(\frac{1}{2}+it,f \otimes \chi)|^2\,\cdot\, (\frac{x}{q^N})^{\beta^\prime}\,\cdot \,\mathcal{K}_{\infty}(\frac{1}{2}+it,0,\beta^\prime,\chi)\,\cdot\,\min\{1,\frac{1}{S|\log (\frac{x}{q^N})|}\}\,dt$$\\

\noindent
{\bf Theorem 6.1.}\\
$$\lim_{S \to \infty} E(x,S)=0\,\,\,\,\,\,(\mbox{for}\,\,\, x>0)$$

\begin{proof}  
We first show that
$$\lim_{S\to\infty}\int_{\delta_0+iS}^{\beta^\prime+iS}\frac{Z(w^\prime)\,x^{w^\prime}}{w^\prime}\,dw^\prime=0\,\,\,\,\,\,\mbox{and}\,\,\,\,\,\,\lim_{S\to\infty}\int_{\delta_0-iS}^{\beta^\prime-iS}\frac{Z(w^\prime)\,x^{w^\prime}}{w^\prime}\,dw^\prime=0$$
Let $w^\prime=\delta+iS$. Then 
$$Z(w^\prime)\ll S^m,\,\,\,\,\,\,(\mbox{for}\,\,\,m<\frac{1}{2}\,\,\,\,\mbox{and}\,\,\,\,|w^\prime|=\sqrt{\delta^2+S^2}\ll S)$$
Thus the integrals above approach $0$ as $S \to \infty$.\\

\noindent
Consider the sets:\\
$A=\Bigl\{N: \frac{1}{S|\log(\frac{x}{q^N})|} \leq \frac{1}{\sqrt{S}}\Bigr\}$\\
$B=\Bigl\{N: \frac{1}{S|\log(\frac{x}{q^N})|} \geq \frac{1}{\sqrt{S}}\Bigr\}$\\

\noindent
On $A$,
$$E(x,S) \ll \frac{1}{\sqrt{S}}\,\sum_{\chi}\,\int_{-\infty}^\infty |L(\frac{1}{2}+it,f \otimes \chi)|^2\,\cdot\, (\frac{x}{q^N})^{\beta^\prime}\,\cdot \,\mathcal{K}_{\infty}(\frac{1}{2}+it,0,\beta^\prime,\chi)\,dt$$
$$=\frac{x^{\beta^\prime}}{\sqrt{S}}\,\sum_{\chi}\,\int_{-\infty}^\infty |L(\frac{1}{2}+it,f \otimes \chi)|^2\,\cdot\, (q^N)^{-\beta^\prime}\,\cdot \,\mathcal{K}_{\infty}(\frac{1}{2}+it,0,\beta^\prime,\chi)\,dt$$
$$=\frac{x^{\beta^\prime}}{\sqrt{S}}\,Z(\beta^\prime)\,\,\,\,\mbox{where}\,\,\,\,Z(\beta^\prime)\,\,\,\,\mbox{is independent of}\,\,\,\,S$$
So 
$$\lim_{S\to\infty} E(x,S)=0$$
On $B$, $\mathcal{K}_{\infty}(\frac{1}{2}+it,0,\beta^\prime,\chi)$ can be estimated by the analytic conductor:
$$Q(\chi,t)=\prod_{v \approx \mathbb{R}}(1+|t+t_v|)\,\cdot\, \prod_{v \approx \mathbb{C}}(1+\ell_v^2+4(t+t_v)^2)$$
Break up $E(x,S)$ into two sums over $q^N\leq \log S$ and $q^N \geq \log S$. Since $Z(w^\prime)$ converges absolutely for $\Re(w^\prime)>1$, the second sum over $q^N \geq \log S$ approaches $0$. So consider
$$\sum_{\chi: q^N \leq \log S}\,\int_{-\infty}^\infty |L(\frac{1}{2}+it,f \otimes \chi)|^2\,\cdot\, (\frac{x}{q^N})^{\beta^\prime}\,\cdot \,\mathcal{K}_{\infty}(\frac{1}{2}+it,0,\beta^\prime,\chi)\,\cdot\,\min\{1,\frac{1}{S|\log (\frac{x}{q^N})|}\}\,dt$$
Now in B,
$$\sum_{q^N \leq \log S}\,\int_{-\infty}^\infty1 \ll (\log S)^k,\,\,\,k>0$$
The convexity bound in the depth aspect gives
$$L(\frac{1}{2}+it,f \otimes \chi) \ll (q^N)^{\frac{1}{2}} \leq (\log S)^{\frac{1}{2}}$$
Fix $\chi=1$ and $0<t<1$ for $v|\infty$. Then
$$\mathcal{K}_{\infty}(\frac{1}{2}+it,0,\beta^\prime,\chi) \ll 1$$
Also
$$\frac{1}{S|\log(\frac{x}{q^N})|} \geq \frac{1}{\sqrt{S}} \Longrightarrow \frac{1}{\sqrt{S}|\log (\frac{x}{q^N})|} \geq 1 \Longrightarrow xe^{-\frac{1}{\sqrt{S}}} \leq q^N \leq xe^{\frac{1}{\sqrt{S}}}$$
This restricts $N$ to a set of measure $\ll \frac{1}{\sqrt{S}}$. So in the second case
$$\lim_{S\to\infty} E(x,S)=0$$
\end{proof}

\noindent
By Cauchy's theorem,
$$\frac{1}{2\pi i}\int_R \frac{Z(w^\prime)\,x^{w^\prime}}{w^\prime}\,dw^\prime=xP(\log x)$$
Indeed, $Z(w^\prime)$ has a pole of order $1$ at $w^\prime=1$, so by the residue theorem:\\
\newcommand{\Res}{Res}
$$\frac{1}{2\pi i}\int_R \frac{Z(w^\prime)\,x^{w^\prime}}{w^\prime}\,dw^\prime=\Res_{w^\prime=1}\Bigl(Z(w^\prime)\,\cdot\,\frac{x^{w^\prime}}{w^\prime}\Bigr)$$
Suppose the Laurent expansion of $Z(w^\prime)$ 
$$Z(w^\prime)=\sum_{n=-\infty}^{\infty} a_n\,(w^\prime -1)^n$$
and
$$x^{w^\prime}=xe^{(w^\prime-1)\log x}=x\,\sum_{n=0}^\infty \frac{(w^\prime-1)^n\,\log^n x}{n!}$$
Then the coefficient of $(w^\prime-1)^{-1}$ in the product $Z(w^\prime)\,\cdot\,\frac{x^{w^\prime}}{w^\prime}$ is $xP(\log x)$, where $P(\log x)$ is a polynomial in $\log x$. So
$$\frac{1}{2\pi i}\int_R \frac{Z(w^\prime)\,x^{w^\prime}}{w^\prime}\,dw^\prime=\frac{1}{2\pi i}\int_{\beta^\prime-iS}^{\beta^\prime+iS} \frac{Z(w^\prime)\,x^{w^\prime}}{w^\prime}\,dw^\prime\,-\,\frac{1}{2\pi i}\int_{\delta_0-iS}^{\delta_0+iS} \frac{Z(w^\prime)\,x^{w^\prime}}{w^\prime}\,dw^\prime$$
$$=xP(\log x)$$
Now Perron's formula showed that
$$\frac{1}{2\pi i}\int_{\beta^\prime-i\infty}^{\beta^\prime+i\infty} \frac{Z(w^\prime)\,x^{w^\prime}}{w^\prime}\,dw^\prime= \sum_{q^N\leq x}\,\int_{-\infty}^\infty |L(\frac{1}{2}+it,f \otimes \chi)|^2\,\cdot\,\mathcal{K}_{\infty}(\frac{1}{2}+it,0,\beta^\prime,\chi)\,dt$$
Thus as $S \to \infty$,
$$\sum_{q^N \leq x}\,\int_{-\infty}^\infty |L(\frac{1}{2}+it,f \otimes \chi)|^2\,\cdot\,\mathcal{K}_{\infty}(\frac{1}{2}+it,0,\beta^\prime,\chi)\,dt=xP(\log x)\,+\,\frac{1}{2\pi i}\int_{\delta_0-i\infty}^{\delta_0+i\infty} \frac{Z(w^\prime)\,x^{w^\prime}}{w^\prime}\,dw^\prime$$\\

\noindent
{\bf Theorem 6.2.}\\
$$\frac{1}{2\pi i}\int_{\delta_0-i\infty}^{\delta_0+i\infty}\frac{Z(w^\prime)\,x^{w^\prime}}{w^\prime}\,dw^\prime \ll x^{\frac{2\delta_0+1}{3}}\,\cdot\, \log x\,\,\,\,\,\,\,\,\,\,\, (\frac{11}{18}<\delta_0<1)$$

\begin{proof}
By the choice of $\delta_0$,
$$\frac{Z(w^\prime)}{w^\prime}=\frac{Z(\delta_0+i\eta)}{\delta_0+i\eta}$$
is a square integrable function on $\mathbb{R}$. Let
$$E(x)=\frac{1}{2\pi i}\int_{\delta_0-i\infty}^{\delta_0+i\infty}\frac{Z(w^\prime)\,x^{w^\prime}}{w^\prime}\,dw^\prime$$\\

\noindent
{\bf Lemma 6.3}\\
$$\int_0^x |E(t)|^2\,dt \ll x^{2\delta_0+1}$$

\begin{proof}
Let $x=e^{-2\pi u}$ and again $w^\prime=\delta+i\eta$. So
\begin{align*}
E(e^{-2\pi u})&=\frac{1}{2\pi i}\int_{-\infty}^{\infty}\frac{Z(\delta_0+i\eta)}{\delta_0+i\eta}\,\cdot\,e^{-2\pi u(\delta_0+i\eta)}\,\cdot\,i\,d\eta\\
&= \frac{1}{2\pi}\int_{-\infty}^{\infty} e^{-2\pi iu\eta}\,\cdot\, f(\eta)\,\cdot\,e^{-2\pi u\delta_0}\,d\eta\,\,\,\,\,\,\,\mbox{(where}\,\,\,\,\,\,f(\eta)=\frac{Z(\delta_0+i\eta)}{\delta_0+i\eta})
\end{align*}
Now
$$\hat{f}(u)=\int_{-\infty}^\infty f(\eta)\,e^{-2\pi i \eta u}\,d\eta$$
Thus
$$e^{2\pi u \delta_0}\,\cdot\,2\pi\,\cdot\,E(e^{-2\pi u})=\hat{f}(u)$$
Using Plancherel's theorem:
$$\int_{-\infty}^\infty |\hat{f}(u)|^2\, du = \int_{-\infty}^\infty |f(\eta)|^2\, d\eta \ll 1$$
So
$$1 \gg 4\pi^2 \int_{-\infty}^\infty |e^{2\pi u \delta_0}\,\cdot\,E(e^{-2\pi u})|^2\, du$$
Replace $e^{-2\pi u}$ by $y$ to get\\

$1 \gg \frac{4\pi^2}{2\pi} \int_{0}^\infty y^{-2 \delta_0}\,\cdot\,|E(y)|^2\,\frac{dy}{y}\,=\,2\pi \int_{0}^\infty y^{-(2 \delta_0+1)}\,\cdot\,|E(y)|^2\,dy$\\

$\geq \int_0^x y^{-(2 \delta_0+1)}\,\cdot\,|E(y)|^2\,dy\,\geq x^{-(2\delta_0+1)}\int_0^x |E(y)|^2\,dy\,\,\,\,\,\,\,\,\mbox{for}\,\,\,0\leq y \leq x$\\

Thus
$$\int_0^x |E(y)|^2\,dy \ll x^{2\delta_0+1},\,\,\,\,\,\,\,0<\delta_0<1$$

\end{proof}

\noindent
We now prove Theorem 6.2, that $$E(x)\ll x^{\frac{2\delta_0+1}{3}}\,\cdot\, \log x$$
First note that $\mathcal{K}_{\infty}(\frac{1}{2}+it,0,\beta^\prime,\chi)$ is positive. Now for $x \leq y$
$$\{N: q^N \leq x\} \subseteq  \{N: q^N \leq y\}$$
Again
$$E(x)=\sum_{q^N \leq x}\,\int_{-\infty}^\infty |L(\frac{1}{2}+it,f \otimes \chi)|^2\,\cdot\,\mathcal{K}_{\infty}(\frac{1}{2}+it,0,\beta^\prime,\chi)\,dt\,-\,xP(\log x)$$
So\\

$E(y)-E(x)=\displaystyle \sum_{q^N \leq y}\,\displaystyle \int_{-\infty}^\infty |L(\frac{1}{2}+it,f \otimes \chi)|^2\,\cdot\,\mathcal{K}_{\infty}(\frac{1}{2}+it,0,\beta^\prime,\chi)\,dt\, - \,$\\
$$\sum_{q^N \leq x}\,\int_{-\infty}^\infty |L(\frac{1}{2}+it,f \otimes \chi)|^2\,\cdot\,\mathcal{K}_{\infty}(\frac{1}{2}+it,0,\beta^\prime,\chi)\,dt\,- \,(yP(\log y)\,-\,xP(\log x))$$
Since $\mathcal{K}_{\infty}(\frac{1}{2}+it,0,\beta^\prime,\chi)$ is positive, 
$$E(y)-E(x) \geq -(yP(\log y) - xP(\log x))$$
Fix $x \geq 3$\\
(a) Replace $y$ with $x+u$ for $0 \leq u \leq x$:\\
$$E(x+u)-E(x) \geq -[(x+u)P\log(x+u) - xP(\log x)]$$
$$\Longrightarrow E(x) \leq  E(x+u)+(x+u)P\log(x+u) - xP(\log x)$$
Now $P$ is a linear polynomial, so rewrite\\

$(x+u)P\log(x+u) - xP(\log x)\,=\,(x+u)[A\log(x+u)+B] - x[A\log x+B]$\\

$=Ax(\log\frac{x+u}{x}) + Au\log(x+u)+Bu\,=\,Ax(\log(1+\frac{u}{x})) + Au\log(x+u)+Bu$\\

Now\\

$\log(1+h) \leq h \,\,\, \mbox{for} \,\,\, 0 \leq h \leq 2\,\Longrightarrow \log(1+\frac{u}{x}) \leq \frac{u}{x} \leq 1$\\

So\\

$Ax(\log(1+\frac{u}{x})) + Au\log(x+u)+Bu\,\leq Ax \cdot \frac{u}{x} + Au\log(x+u)+Bu$\\

$=Au+Bu+Au\log(x+u)\,\leq Du\log x + Au(\log x+\log 2)\,\,\,\mbox{since}\,\,\,u \leq x$\\

Thus
$$E(x) \leq E(x+u) +Cu \log x \,\,\, \mbox{for some constant}\,\,\,C$$

\noindent
(b) Replace $x$ with $x-u$ and $y$ with $x$ for $0 \leq u<x$. Then
$$E(x)-E(x-u) \geq -[xP(\log x) - (x-u)P\log(x-u)]$$
$$\Longrightarrow E(x) \geq E(x-u) - Cu \log x$$
\noindent
Let $0 \leq H \leq x$. Integrate the inequalities over $0 \leq u \leq H$:
$$\int_0^H E(x)\,du \leq \int_0^H (E(x+u)+Cu \log x) \, du\,=\, H \cdot E(x) \leq \int_0^H E(x+u)\,du + \frac{C}{2}\,H^2\,\log x$$
and
$$H \cdot E(x) \geq \int_0^H E(x-u)\,du - \frac{C}{2}\,H^2\,\log x$$
So
$$\int_0^H E(x-u)\,du - \frac{C}{2}\,H^2\,\log x \leq H \cdot E(x) \leq \int_0^H E(x+u)\,du + \frac{C}{2}\,H^2\,\log x$$
Change variables and replace $\frac{C}{2}$ with $C$ to get
$$\frac{1}{H}\, \int_{x-H}^x E(t)\, dt\,-\, CH \log x \leq E(x) \leq \frac{1}{H}\, \int_x^{x+H} E(t)\, dt\,+\, CH \log x $$
For $E(x) \geq 0$, apply the second inequality, otherwise apply the first. So for $E(x) \geq 0$, 
$$E(x)^2 \ll  \frac{1}{H^2}\, \Bigl(\int_x^{x+H} E(t)\, dt \Bigr)^2\,+\, C^2H^2 \log^2 x $$
Apply Cauchy-Schwarz:
$$E(x)^2 \ll  \frac{1}{H^2}\,\int_x^{x+H} |E(t)|^2\, dt \,\cdot \int_x^{x+H} 1\,dt \,+\, H^2\log^2 x $$
$$=\frac{1}{H}\,\int_x^{x+H} |E(t)|^2\, dt \,+\, H^2\log^2 x \ll \frac{1}{H}\,\cdot\, x^{2\delta_0+1}\,+\,H^2 \log x$$
since
$$\int_0^x |E(t)|^2\,dt \ll x^{2\delta_0+1}\,\,\,\mbox{and}\,\,\,H\leq x$$
We want $\frac{1}{H}\,\cdot\, x^{2\delta_0+1}=H^2$, so take $$H=x^{\frac{2\delta_0+1}{3}}$$ Then
$$E(x) \ll H \log x = x^{\frac{2\delta_0+1}{3}}\,\cdot\, \log x$$
\end{proof}

\noindent
Recall that the $\chi$-depth-aspect convexity bound for the twisted $L$-function is
$$L(\frac{1}{2}+it,f \otimes \chi) \ll q^{N(\frac{d}{2}+\epsilon)}$$
Let us now use the results obtained above to break convexity by decreasing the exponent of $q^N$.\\

\noindent
Choose $H$ such that $$x^{\frac{2\delta_0+1}{3}} \ll H \ll x^{\frac{2\delta_0+1}{3}}$$
Let
 
$S(x)=\displaystyle \sum_{q^N \leq x}\,\displaystyle \int_{-\infty}^\infty |L(\frac{1}{2}+it,f \otimes \chi)|^2\,\cdot\,\mathcal{K}_{\infty}(\frac{1}{2}+it,0,\beta^\prime,\chi)\,dt$\\

$=xP(\log x)\,+\,O(x^{\frac{2\delta_0+1}{3}}\,\log x)\,=\,x P(\log x)\,+\,E(x)$\\

\noindent
Now for $H>0$, $\{N:q^N \leq x \} \subset \{N:q^N \leq x+H \}$ and \mbox{$\mathcal{K}_{\infty}(\frac{1}{2}+it,0,\beta^\prime,\chi)$} is {\it positive}. So for trivial $\chi$,
$$S(x+H+1)-S(x) \geq \sum_{x \leq q^N \leq x+H}\int_{-\infty}^\infty |L(\frac{1}{2}+it,f \otimes \chi)|^2\,\cdot\,\prod_{v|\infty}\mathcal{K}_v(\frac{1}{2}+it,0,\beta^\prime,1)\,dt$$ 
Now
$$S(x+H+1)-S(x)=(x+H+1)P(\log(x+H+1))-xP(\log x)+ E(x+H+1)-E(x)$$
Since $x^{\frac{2\delta_0+1}{3}} \ll H \ll x^{\frac{2\delta_0+1}{3}}$ and $E(x) \ll x^{\frac{2\delta_0+1}{3}}$ then
$$E(x+H+1)-E(x) \ll x^{\frac{2\delta_0+1}{3}}\,\log x$$
and
$$(x+H+1)P(\log(x+H+1))-xP(\log x) \leq C(H+1) \log x$$
So
$$S(x+H+1)-S(x) \ll x^{\frac{2\delta_0+1}{3}}\,\cdot \,\log x$$
$$\Longrightarrow \sum_{x \leq q^N \leq x+H} \int_{-\infty}^\infty |L(\frac{1}{2}+it,f \otimes \chi)|^2\,\cdot\,\prod_{v|\infty}\mathcal{K}_v(\frac{1}{2}+it,0,\beta^\prime,1)\,dt \ll x^{\frac{2\delta_0+1}{3}}\,\cdot\, \log x$$
Now
$$Q(\chi,t)^{-\beta^\prime} \ll \mathcal{K}_v(\frac{1}{2}+it,0,w,\chi) \ll Q(\chi,t)^{-\beta^\prime}\,\,\,\,\,\,(\mbox{for}\,\,\,v|\infty)$$
where
$$Q(\chi,t)=\prod_{v \approx \mathbb{R}}(1+|t+t_v|)\,\cdot\, \prod_{v \approx \mathbb{C}}(1+\ell_v^2+4(t+t_v)^2)$$
For trivial $\chi$, $t_v=l_v=0$. Also recall for $v|\infty$, fix $0<t<1$. Then
$$\mathcal{K}_v(\frac{1}{2}+it,0,\beta^\prime,1) \gg (1)^{-(d-1)\beta^\prime}$$
So
$$x^{\frac{2\delta_0+1}{3}}\,\log x \gg \sum_{x \leq q^N \leq x+H}\int_{-\infty}^\infty |L(\frac{1}{2}+it,f \otimes \chi)|^2\,\cdot\,\prod_{v|\infty}\mathcal{K}_v(\frac{1}{2}+it,0,\beta^\prime,1)\,dt$$
$$\gg \sum_{x \leq q^N \leq x+H} |L(\frac{1}{2}+it,f \otimes \chi)|^2\,\cdot\,(1)^{-(d-1)\beta^\prime}\,dt$$
$$\gg \sum_{x \leq q^N \leq x+H}  |L(\frac{1}{2}+it,f \otimes \chi)|^2\,\cdot\,(q^N)^{-(d-1)\beta^\prime}\,dt$$
$$\geq \sum_{x \leq q^N \leq x+H} |L(\frac{1}{2}+it,f \otimes \chi)|^2\,\cdot\,(x+H)^{-(d-1)\beta^\prime}$$
So
$$\sum_{x \leq q^N \leq x+H} |L(\frac{1}{2}+it,f \otimes \chi)|^2\,dt \ll (x+H)^{(d-1)\beta^\prime}\,\cdot \, x^{\frac{2\delta_0+1}{3}}\,\cdot \,\log x$$
$$\ll x^{d-1+\frac{2\delta_0+1}{3}+\frac{\epsilon}{2}}\,\cdot \,\log x\,\,\,\mbox{where}\,\,\,\beta^\prime=1+\frac{\epsilon}{2d-2}$$
$$\ll x^{d-1+\frac{2\delta_0+1}{3}+\epsilon}$$
By a standard argument analogous to that in [Good 1982], this short-interval moment bound implies the pointwise bound
$$L(\frac{1}{2}+it,f \otimes \chi)\ll (q^N)^{\frac{d-1}{2}+\frac{2\delta_0+1}{6}+\epsilon} \ll (q^N)^{\frac{d-1+\vartheta}{2}+\epsilon}$$
for $\vartheta <1$.

\end{document}